%% file: n111.ltx
\documentclass[12pt,a4paper,twoside]{article}

\usepackage{varioref}

\usepackage[dvips]{graphicx} 

\usepackage{color}

\usepackage[hyperindex]{hyperref}

\usepackage[T1,OT1]{fontenc}

\usepackage[latin1]{inputenc}

 \usepackage{fancyhdr}

 \fancypagestyle{plain}{%
 \fancyhf{}%
 \fancyhead[LE,RO]{\bfseries\thepage}}

 \fancypagestyle{index}{%
 \fancyhf{}%
 \fancyhead[LE,RO]{\bfseries\thepage}
 \fancyhead[LO]{\bfseries Index}
 \fancyhead[RE]{\bfseries Index}}

\usepackage{makeidx}

\usepackage{float}

\usepackage[numindex,nottoc,numbib,section]{tocbibind} 
 
\usepackage[titles]{tocloft}    

\usepackage{enumerate}

\usepackage{array}

\usepackage[size=2em]{diagrams} 

\usepackage[leqno]{amsmath}     
\usepackage{amsthm}             
\usepackage{amssymb}            
\usepackage{mathrsfs}           
\usepackage{dsfont}             
\DeclareSymbolFontAlphabet{\mathbbb}{AMSb}

\newtheorem{theorem}{Theorem}
\newtheorem{lemma}[theorem]{Lemma}
\newtheorem{proposition}[theorem]{Proposition}
\newtheorem{definition}[theorem]{Definition}

\newtheorem{cor}[theorem]{Corollary}

\theoremstyle{definition}

\newtheorem*{remark}{Remark}

\newcommand{\bs}{\boldsymbol}

\newcommand{\abs}[1]{\lvert#1\lvert}

\makeindex

\title{Ramunajan $(n_1,n_2,\ldots,n_{d-1})$-regular hypergraphs based on\\ 
  Bruhat-Tits Buildings of type $\tilde{A}_{d-1}$} 

\author{Alireza Sarveniazi}

\date{02.01.2004} 
\begin{document}
\maketitle
\ 
\vfill
\
\tableofcontents \newpage


\pagestyle{plain} 

\input{n112.ltx}

 \clearpage 

\bibliographystyle{plain}

\input{n111.bbl}
\nocite*                    
\clearpage 

\end{document}

%% file: n112.ltx
\section{Introduction} 
               
  We introduce in this text Ramanujan regular hypergraphs, generalizing the
  definition of Ramanujan graphs introduced by Lubotzky, Philips, Sarnak
  \cite{LPS1}.  
   Ramanujan graphs are regular graphs, whose adjacency matrices,
  $($or equivalently their Laplacians$)$, have eigenvalues satisfying some
  natural upper bounds. Graphs satisfying such bounds have many interesting
  properties.  For example they are very powerful expander graphs, which make
  them interesting objects for many applications from communication networks to
  computer science. The interesting situation here, is that on the one hand it
  is known and easy to see that almost all regular graphs in a precise sense 
  are Ramanujan graphs. On the other hand it is difficult to show for an
  explicitly given regular graph that it is a Ramanujan graph. It is therefore
  quite surprising and has lot of interest, that some very explicit
  graphs coming from number theory can be shown to be Ramanujan graphs. 
 These graphs are constructed using unit groups of quaternion algebra and
  their action on the associated symmetric spaces, which in these examples
  just are the trees in the sense of \cite{Serr1} or \cite{tits2}.  The quotients
  of these trees by the action of the unit groups mentioned above give
  already  the Ramanujan graphs we are discussing. To show the Ramanujan property
  requires deep tools from the theory of automorphic representations. In fact
  it turns out that the Ramanujan property of these graphs is equivalent to
  the fact, that the associated representations satisfy the Ramanujan-Petersson
  conjecture. This is not very difficult to prove because one sees rather
  directly that the adjacency operators of the quotient graph above is more or
  less nothing else than the Hecke operator of the corresponding prime.
The first examples given where working with quaternion algebras over the
  rational numbers and were using Deligne's proof of the Ramanujan-Petersson
  conjecture plus the Jaquet-Langlands correspondence between automorphic
  representations of quaternion algebras and cuspidal automorphic
  representations of $GL(2,\mathbb{Q})$ of an appropriate type. A later
  variant of this by M.Morgenstern \cite{mor1} was using instead quaternion
  algebras over the rational function field $\mathbb{F}_q(t)$ and
  corresponding results of Drinfeld \cite{Drin1} for the associated automorphic representations. It is
  this class of examples of Morgenstern, which we will generalize to higher
  dimensions. This is made possible by the observation, that the
  quaternion algebras in Morgenstern's examples are nothing else than the
  quotient $($skew$)$ fields of skew polynomial rings which are well known
  in the theory of Drinfeld modules. So we consider the skew polynomial ring
  $\mathbb{F}_{q^d}\{\tau\}$ over the field $\mathbb{F}_{q^d}$ of $q^d$ elements,
  where the indeterminate $\tau$ satisfies the rule $\tau \cdot \lambda =
  \lambda^q \cdot\tau$ for $\lambda \in \mathbb{F}_{q^d}$.\\ The center of this
  ring is the polynomial ring $\mathbb{F}_q[t]$, the $($skew$)$ quotient field
  $\mathbb{F}_q(\tau)$ is a division algebra of dimension $d^2$ over the
  center $\mathbb{F}_q(t)$, the field of rational functions over
  $\mathbb{F}_q$.\\ This division algebra is ramified exactly at the primes
  $t=0$ and $t=\infty$. We consider the localization
  $\mathbb{F}_{q^d}\{\tau\}[\frac{1}{p}]$ at a prime $p = p(t)$ of
  $\mathbb{F}_q[t]$ different from $t$. The unit group associated with this
  maximal order over the center
  $\mathbb{F}_q[t][\frac{1}{p}]$ is acting on the Bruhat-Tits building at the
  prime $p$, which is a building of type $\tilde{A}_{d-1}$. The quotients by
  the unit groups we define, will be simplicial complexes of dimension
  $(d-1)$. It is these simplicial complexes and their adjacency matrices, we
  are studying here. In this text, we partly work in the category of
  simplicial complexes and more specialized buildings and partly in the
  category of hypergraphs. The procedure to be followed now, is  quite similar
  to the case of graphs. It consists in identifying the adjacency operator of
  these quotient complexes with some Hecke operator acting on certain 
  corresponding spaces of automorphic forms and then making use of the relevant
  results about automorphic representations in this context. The main result
  we are relying on here, is Lafforgue's recent proof of the Langlands
  conjecture for $GL(n)$ over a function field and in particular the
  Ramanujan-Petersson conjecture for these cuspidal automorphic
  forms. What causes difficulties, is that here the Jacquet-Langlands
  correspondence between automorphic representations for  $GL(n)$ and
  automorphic representations for the division algebra is not effective
  enough. Using this correspondence, it is possible that noncuspidal
  automorphic representations come up. Lafforgue's result however requires the
  automorphic representation to be cuspidal.\\
 Nevertheless, using additionally  \cite[Weak
  Lifting theorem, $4.2$]{Artur} it is
  possible to concude the Ramanujan property at least for the d a prime from
  Lafforgue's result. \\
It should be metioned also, that after my Ph.D.thesis was $($more or less$)$
  finished, preprints of W.-C. Winnie Li \cite{wlij} and Lubotzky, Samuels and
  Vishne \cite{LSV2},  \cite{LSV1}  become avilable to us. In particular 
  W. Lie  using a trick of L.Clozel, reduce the the
  Ramanujan property of the automorphic representation to the  Ramanujan
  property of the moduli scheme of $\mathcal{D}$ elliptic sheaves \cite{LRS}.

Let us add at this point, that there is another treatment of higher dimensional
  situations in the literature, namely \cite{BL1}  B.W.Jordan and R.Livne. These
  results work again with unit groups of quaternion algebras, but allow
  denominators at more than one prime. Correspondingly, their unit groups are
  acting on a product of trees. As this is again the $\text{\rm{GL}}(2)$-case, the
  Ramanujan-Petersson property holds. Also C.Ballantine \cite{Ba} has given a
  hypergraph for $ d= 3$ related to the Bruhat-Tits building associted to $\text{PGL}(3; \mathbb{Q}_p)$. \\
We end this introduction with a
  description of the content of the different chapters of this text.\\

Section $2$ describes buildings of type $\tilde{A}_n$  and various concepts
  related to this. A building is a geometrical-combinatorial object.  
The definition of affine building is strictly combinatorial, understanding the relation to matrix groups is one of the
keys to understanding affine buildings.  In this section we discuss the affine buildings associated to
the matrix groups PGL(d,F), where $F$ is a local field.\\ The aim of this
  section  is to give an introduction and as elementary as possible to the subject of buildings of
  type $\tilde{A}_{d-1}$ related to non-archimeadian local fields.

In section $3$ we explain our situation and basic definitions. $4.$ gives  a short review of
some fundamental facts concerning automorphic representations. In particular
the Hecke algebra and the standard generators of the local Hecke algebras will
be explained. $5$ explains identification between  the adjacency operators of the quotient complex
resp. quotient hypergraph with the corresponding Hecke operators of the
automorphic situation. Finally the Ramanujan property is discussed in the section$6$. In section $6$ we describe our examples in general. We
  indicate some of of the concepts from the theory of automorphic
  representations. In particular we describe the relationship between
  adjacency operators and Hecke operators. Finally we discuss in this chapter,
  what would be needed additionally from the theory of automorphic
  representations, in particular regarding the work of L.Lafforgue.
In this section we describe our candidates for Ramanujan hypergraphs in
general. The main point to be explained is the relation between the
different adjacency operators in a hypergraph and the corresponding Hecke
operators for certain automorphic representations, which are related. To
really conclude the Ramanujan property, one has to use the relevant results
for automorphic representations. Here the case $d=2$ is easier than the case
of general $d$, because only in this case there is a fully worked out
Jacquet-Langlands correspondence between automorphic representations of
unit groups of division algebras and  general linear groups over, which should allow to relate the situation to automorphic
representations of the general linear group and thereby to the recent results
of Lafforgue, Instead we use a result of W.Li resp. L. Clozel, which work with
the  theory of $\mathcal{D}$-elliptic sheaves and does not make use of the
Jacquet-Langlands correspondence.\\ In section $7$ we consider the main properties of skew polynomial rings
$\mathbb{F}_{q^d}\{\tau\}$ and some related rings. Section $8$ describes the
arithmetic groups we want to study. Sections $9$ and $10$ gives the explicit
construction of the Ramanujan hypergraphs, our simplicial complexes we are
introduced in. Again these will be described in terms of Cayley graphs of
various groups.   We will giveiIn section   $10$ a very simple but
combinatorical description of our hypergraph as Cayley graph of group
$\text{PGL}(d, F_0)$ or $\text{PSL}(d, F_0)$ over a finite field $F_0$.
 \paragraph{Aknowledgement} \
\\
 I am indebted to Prof. Dr. Ulrich Stuhler, who supervised my doctoral work,
 for helpful discussion, continued   encouragement and suggestions for
 presentation in the form my thesis and this text. I am also much thankful to
 PD Dr. Thomas Lehmkuhl for useful discussion, specially on Theorem
 \ref{above108}. I am also very thankful to Prof. Dr. Larry Smith for his
 careful reading of this text. 
 \section{Affine Buildings}\label{2Baffine}
In this section $F$ denotes a local non archimedean field with valuation ring
$\mathcal{O}$, $(\pi)$ is the maximal ideal with specified uniformizing element
$\pi \in \mathcal{O}$.\\ $F^d$ denotes the $d$-dimensional  vector
space over $F$. 
\begin{definition}
$i)$ A lattice in $F^d$ is a $\mathcal{O}$-submodule of 
$F^d$, for which  there exists a basis $\{b_1,\ldots,b_d\}$ of $F^d$ such that $
L=\sum_{i=1}^d\mathcal{O}b_i$ holds.\\
$(ii)$ Two lattices $L,L^\prime \subset F^d$ are similar iff there exists
$\lambda \in F^\times$, such that $L^\prime=\lambda L$. In this situation we
denote $L\sim L^\prime$.
\end{definition}
\begin{remark}
$(i)$ If L is a finitely generated torsion free $\mathcal{O}$-submodule of 
$F^d$ satisfying $F\cdot L=F^d$, then $L$ is a lattice in the sense of the definition
above.\\$(ii)$ Of course a similar definition as above can be made for arbitrary
abstract vector spaces $V$ over $F$ of finite dimension.
\end{remark} The group $\text{\rm{GL}}(d, F)$ of $F$-linear isomorphisms of the vector
space $F^d$ acts on the set of lattices by \[ g L=\{g(l) | l\in L\}\]
for $g\in \text{\rm{GL}}(d, F)=Aut_F(F^d), L\subset F^d$ an
$\mathcal{O}$-lattice.\begin{lemma} The action of $\text{\rm{GL}}(d, F)$ is transitive on
  $\mathcal{L}$ where $\mathcal{L}$ is the set of all lattices.\end{lemma}
\begin{proof} Denote $L_0:= \bigoplus_{i=1}^d \mathcal{O}e_i$ the standard lattice
  $\mathcal{O}^d\subset F^d,\{e_1,\ldots,e_d\}$ the standard basis of $F^d$.\\
It is sufficient to show, that for an arbitrary lattice $L \in \mathcal{L}$
there exists $g\in \text{\rm{GL}}(d, F), g.L_0= L$.\\ By definition, there exists a basis
$\{b_1,\ldots,b_d\}$ of $F^d$, such that
$L=\sum_{i=1}^d\mathcal{O}b_i$. Define $g\in \text{\rm{GL}}(d, F)$ by $g(e_i)=b_i$ for
$i=1,\ldots,d$. Obviously $g.L_0= L$ holds.
\end{proof}
\begin{lemma} The action of $\text{\rm{GL}}(d, F)$ on $\mathcal{L}$ induces a
  corresponding action on $\mathcal{L}/\sim$, which is again transitive.\end{lemma}
\begin{proof} It is enough to remark , that $L\sim L^\prime$ implies $g.L\sim
  g.L^\prime$ for any element $g\in \text{\rm{GL}}(d, F)$.
\end{proof}
\begin{definition} A simplicial complex $X_{\bs{\cdot}}$ is given as\\
$(i)$ a set $X_0$ $($the set of vertices of the simplicial complex $X_{\bs{\cdot}})$\\
$(ii)$  For any natural number $d \geq 0$, the set of $d$-simplexes $X_d$ where
$X_d \subset \mathbb{P}(X_0)$, the power set of  $X_0$  and for $ Y \in X_d$ one
has  $ \abs{ Y } =d+1$ for the cardinality  
$ \abs{Y}$  of $Y$.\\
$(iii)$ If  $ Y \in X_d$ and $ Y^\prime \subset Y$, $Y^\prime \neq \varnothing$, then
  $ Y^\prime \in  X_{d^\prime}, d^\prime \leq d$, where  $ d^\prime + 1 =
  \abs{ Y^\prime }$.
\end{definition}
\begin{remark}
$(a)$ $(iii)$ says the following: \\
 If $Y$ is a $d$-simplex, any nonempty subset $Y^\prime$ of 
$(d^\prime + 1)$-elements is a  $d^\prime$-simplex of $X_.$.\\
$(b)$ There are various concepts related  to the concept of a simplicial complex
 as ordered simplicial complex and simplicial set. To all of these, one can associate
 functorially a topological space $\abs{X_{\bs{\cdot}}}$ $($the so called realization of
 $X_{\bs{\cdot}})$ and  $X_{\bs{\cdot}}$ gives in a certain sense a combinatorial description of the topological space $\abs{X_{\bs{\cdot}}}$.  For all of these materials  see
 \cite[chapter 8]{Ro1}. \\
We define now the so called affine building
 associated with the group $\text{\rm{GL}}(d,F)$ (or if one prefers, the projective
 linear group $\text{PGL}(d,F)$). It will be a simplicial complex of dimension
$(d-1)$.
\end{remark}
 \begin{definition}
 The affine building $X_{\bs{\cdot}} = X_{\bs{\cdot}}(F^d)$ (associated to the group $\text{PGL}(d,F)$) is a
 simplicial complex given as follows:\\
$($i$)$ the set of vertices $X_0$ is $\mathcal{L}/\sim$, the set of lattices up
 to similarity in $F^d$.\\
$($ii$)$ If  \[L_0 \supset L_1 \supset \ldots L_r \supset \pi L_0\]
is a flag of $(r + 1)$ different lattices in  $F^d$, then $< L_0,L_1, \ldots,
 L_r>$ is a $r$-simplex in $X_{\bs{\cdot}}$, that is   $< L_0,L_1, \ldots,
 L_r> \in X_r$. Any simplex of $X_{\bs{\cdot}}$ is obtained in this way.
\end{definition}
\begin{remark}
$($i$)$ As $\dim_k{(L_0/\pi L_0)} = d)$, where $k=\mathcal{O}/\pi\mathcal{O}$ is the
residue field, it follows:\\
 $\dim{(X_{\bs{\cdot}})}=d-1$, that is, the maximum dimension of a simplex in
 $X_.$ is $d-1$.
$($ii$)$ There is an obvious action of $\text{\rm{GL}}(d,F)$ on $X_{\bs{\cdot}}$, given by:
\[ g< L_0,L_1, \ldots,L_r> = < gL_0,gL_1, \ldots,gL_r> \]
for $g\in \text{\rm{GL}}(d,F)$, $< L_0,L_1, \ldots,L_r>\in X_{\bs{\cdot}}$. This
action induces an action of $\text{PGL}(d,F)$, because the center of
$\text{\rm{GL}}(d,F)$ acts trivially on $X_{\bs{\cdot}}$.\\
$($iii$)$ In the language of buildings the $(d-1)$-dimensional simplices are
called chambers of  $X_{\bs{\cdot}}$. They are given as complete flags:
\[ L_0 \supset L_1 \supset \ldots L_{d-1} \supset \pi L_0 \] of length $(d-1)$. As is shown
in \cite{Br}, $X_{\bs{\cdot}}$ is a building in the sense of \cite{Br}. \\In particular, the
apartments of $X_{\bs{\cdot}}(F^d)$ are given as follows: if $F^d = \oplus_{i=1}^dW_i$
is a direct decomposition of $F^d$ into one-dimensional linear subspaces,
 \[ <L> \ \in \underline{A}(W_1,\ldots ,W_d)\quad \text{iff} \quad L =
 \oplus_{i=1}^d(W_i\cap L) \]
for a lattice $L$. $\underline{A}(W_1,\ldots ,W_d)$ is the full sub complex of
$X_{\bs{\cdot}}(F^d)$ generated by these vertices. 
In more concrete terms one can give the following description of the set of 
 $\underline{A}(W_1,\ldots ,W_d)$. Suppose  $W_i = F\omega_i$, $(i=1,\ldots,d)$, such
that $\{\omega_1,\ldots,\omega_d\}$ is a basis of $F^d$. Then $<L>\in
\underline{A}(W_1,\ldots ,W_d)$ iff $L =  \oplus_{i=1}^d\mathcal{O}\pi^n_i\omega_i$ for a
system of appropriate $n_1,\ldots,n_d \in \mathbb{Z}$.\\
simplicial complex  $X_{\bs{\cdot}}(F^d)$ is simplicially contractible.\\
$($ii$)$ The associated  topological space $| X_{\bs{\cdot}}(F^d)|$ is  contractible.\end{remark}
\begin{proof}
$($i$)$ Fix the vertex $<L_0> = <\mathcal{O}^d>$. We will describe a contraction of 
 $X_.(F^d)$ towards the vertex $<\mathcal{O}^d>$.  Suppose  $<L>$ is an
 arbitrary vertex of  $X_.(F^d)$. We can assume (up to change by a scalar
 factor) $L_0 \subseteq L  \subseteq \pi^m L_0$ and   $m$ is maximal with this
 property. Then we define the map
 \[ L \longmapsto T(L) : = L + \pi^{m-1}L_0 \quad \text{if}\quad   1\leq
 m.\]If $m = 0$, we define $T(L) = L_0$, as $L = L_0$ exactly holds. One can
 check the following points easily:\\
$1)$ $T$ induces a well defined map $ \overline{T_0} : X_0(F^d) \longrightarrow
 X_0(F^d)$ on the vertices of  $X_.(F^d)$.\\
$2)$  $ \overline{T_0}$ induces a simplicial morphism 
\[  \overline{T}:  X_.(F^d) \longrightarrow
 X_.(F^d). \]
$3)$ $ \overline{T}$ is homotopic to the identity morphism.\\
For a proof of these facts see \cite{Ro1}. This shows $i)$ in the proposition, $ii)$ is
 an immediate consequence. 
\end{proof}
Fixing a vertex $<L>$ in the building  $X_{\bs{\cdot}}(F^d)$, we consider the link
$lk_{X_{\bs{\cdot}}}(L)$ of the vertex  $<L>$ in  $X_{\bs{\cdot}}(F^d)$. This is the simplicial
complex, given by all simplices $\Delta \in X_{\bs{\cdot}}(F^d)$, such that $\Delta$ does
not contain $<L>$ as a vertex, but $\Delta \cup{L}$ is a simplex in
$X_{\bs{\cdot}}(F^d)$.
\begin{proposition}\label{link0}
 The simplicial complex $lk_{X_{\bs{\cdot}}}(L)$ is isomorphic to the Tits building $X_{\bs{\cdot}}(V)$
 associated to the vector space $ V:= L/\pi L \cong k^d$\end{proposition}
\begin{proof}
 The vertices of the Tits building associated to the vector space  $V \cong k^d$ are given as $<W>$, where $0
 \subsetneqq W \subsetneqq V$ are the proper linear subspaces of $V$.  A $r$-simplex is given
 as $<W_0,W_1,\ldots, W_r>$, where $W_0\supsetneqq W_1\ldots \supsetneqq W_r$ and the
 $<W_j>$ are vertices. It is then immediate to see that the morphism
 \[lk_{X_{\bs{\cdot}}}(L)\longrightarrow X_{\bs{\cdot}}(V)\] 
\[ <L_0,\ldots, L_r> \longmapsto <
 L_0/\pi L, L_1/\pi L, \ldots , L_r/\pi L>\]
$($where $L\supset L_0 \supset \ldots \supset L_r \supset \pi L$ holds$)$ is in
fact an isomorphisms of simplicial complexes.
\end{proof}
\begin{cor}
 i) The Tits building of $V$ resp. $lk_{X_{\bs{\cdot}}}(L)$ is a labeled simplicial complex or building in
 the sense of \cite{Br} by the map $<W>\longmapsto dim(W)$ for $<W> \in
 X_{\bs{\cdot}}(k^d)$. The set of labels is $\{1,\ldots, dim(V) - 1\}$.\\
ii)   $lk_{X_{\bs{\cdot}}}(L)$ is a labeled simplicial complex or even a building.
\begin{proof} i) is clear, ii) follows from i) and the above proposition.
\end{proof} 
\end{cor}
\begin{remark} Of course the labeling in ii) depends on the choice of the
 lattice $L$ or the vertex $<L>$. It does not correspond  to a global labeling
 of the whole building $X_{\bs{\cdot}}(F^d)$.\\
On the other hand there is the possibility to give a labeling to
 $X_{\bs{\cdot}}(F^d)$.This is only $ \text{SL}(d,F)$-invariant, not  $ \text{\rm{GL}}(d,F)$-invariant.It is
 obtained by fixing a Haar measure $\mu$ on the commutative locally compact
 abelian group $F^d$ such that $\mu (\mathcal{O}^d) = 1$. For any lattice $L \in
 \mathcal{L}$ there is a lattice $L^\prime \sim L$, such that
 \[ 1 \geq \mu (L^\prime) \ge \frac{1}{|k|^d} \] holds. This follows immediately from the formula
\[\mu (\pi^m L) =  \frac{1}{|k|^{md}}\mu(L)\] Furthermore , one has $|\mu
 (L^\prime)|= |k|^{-d}, \quad 0\leq d^\prime \leq d-1$. To a neighboring
 vertex $L^\prime$ to $L$ one can
 define a labeling $\lambda:X_{\bs{\cdot}}(F^d) \longrightarrow \{0,1,\ldots,d-1\}$,
 $\lambda (L^\prime): = d^\prime$. 
\end{remark}
\begin{remark}
This global labeling of  $X_.(F^d)$ is different from the local labeling
considered before.\end{remark}
\ \\
Suppose now, $\Delta_1,\Delta_2 \in X_.(F^d)$ are two simplices of the
building.

\begin{definition}
A gallery in the building $X_.(F^d)$ is a sequence of chambers
$[C_0,\ldots,C_r]$, such that $C_i,C_{i+1}$ are neighboring chambers in
$X_.(F^d)$ for $i=0,\ldots, r-1$.\end{definition}
\begin{definition}
The combinatorial distance between  $\Delta_1,\Delta_2 \in X_.(F^d)$ is given
as: 
\begin{equation*}
\min{ \{(r-1)\in \mathbb{N}\cup \{ 0\}\quad | \quad \exists \quad\text{a
    gallery} \quad [C_0,\ldots,C_r],\quad \Delta_1 \subset \overline{C_0},\Delta_2 \subset
\overline{C_r} \} }
\end{equation*}
\end{definition}
\begin{proposition}
$\Delta_1,\Delta_2 \in X_.(F^d)$ are given as above, $\underline{A}$ is any
apartment, such that $\Delta_1,\Delta_2 \in \underline{A}$. Then any
shortest gallery  $[C_0,\ldots,C_r]$ as in the definition above satisfying $C_0, C_r
\in \underline{A}$, is contained completely in $\underline{A}$.
\begin{proof}
see \cite[IV.4.Prop., p.88]{Br}
\end{proof}
\end{proposition}

\paragraph{ \text{Description of a standard appartment:}} \ \\
 We want to give an explicit numerical description of an appartment
 $\underline{A}(W_1,\ldots,W_d)$, where the $W_i \subset F^d$    are
 one-dimensional linear spaces and $\oplus_{i=1}^dW_i =F^d$ is a direct sum
 decomposition. It suffices to assume for this purpose, that $W_i = Fe_i$,
 $\{e_1,\ldots,e_d\}$ the standard basis of $F^d$. A vertex $<L> \in
 \underline{A}(W_1,\ldots,W_d)$ in this case iff \[L = \oplus_{i=1}^d(L\cap
 W_i) = (\mathcal{O}\pi^{n_1}e_1\oplus \ldots \oplus \mathcal{O}\pi^{n_d}e_d)
 .\] Therfore $L$ is given by the vector $v(L) = (n_1, \ldots, n_d) \in \mathbb{Z}^d$, a
 vertex $<L>$ is uniquely given by the residue class of $(n_1, \ldots, n_d) \in
 \mathbb{Z}^d/\mathbb{Z}(1,\ldots,1)$. \\ we are going to give a description
 of the chambers $<L_0,\ldots,L_d>$ in the appartment
 $\underline{A}(W_1,\ldots,W_d)$. Denote 
\begin{align*}
v(L_j) &:= (n_1^{(j)},\ldots, n_d^{(j)})\\ 
v(<L_j>) &:= (n_1^{(j)},\ldots, n_d^{(j)}) \quad \mod{\mathbb{Z}(1,\ldots,1)}
\end{align*}
Fixing the lexiographical order on the abelian group $\mathbb{Z}^d$, we can denote
representatives $L_j \in <L_j>$, such that, upon reordering indices $j \in
\{0,\ldots,d\}$, we can assume that 
\begin{equation*}
v(L_0) < v(L_1) < \ldots < v(L_d) < v(L_0) + (1,\ldots,1)
\end{equation*}
holds. This is equivalent to
 \begin{equation*}
L_0 \supset L_1 \supset \ldots < L_d \supset \pi L_0.
\end{equation*}
Finally, if we went to fix types, such that $\mu(L_0) = \mu(\mathcal{O}^d)$,
the standard lattice, we have a unique choice of $j\in \{0,\ldots,d\}$, such
that \begin{equation*}
\log_q{\mu(L_j)} \equiv 0 \quad \mod{(d)}
\end{equation*}
Upon replacing $L_j$ by an appropriate representative $\pi^\alpha L_j$, we can
assume, that \[\mu(\pi^\alpha L_j) = \mu(\mathcal{O}^d) = 1.\] We obtain the
simplices \begin{equation*}
<\pi^\alpha L_j, \pi^\alpha L_{j+1},\ldots, \pi^\alpha L_d, \pi^{\alpha+1}
L_0,\ldots,\pi^{\alpha+1}L_{j-1},\ldots)
\end{equation*}  satisfying 
\begin{align*}
&i) \quad \pi^\alpha L_j \supset \pi^\alpha L_{j+1}\subset \ldots \supset
\pi^{\alpha+1}L_{j-1} \supset \pi^{\alpha+1}L_j\\
&ii) \mu(\pi^\alpha L_j) = 1
\end{align*}
and correspondingly for the other volumes.
\section{Combinatorics and definition of Ramanujan Hypergraph}

\paragraph{Introduction }\ \\
For a finite regular graph, the eigenvalue $\lambda$ of the adjacency matrix
which has the second largest absolute value is of particular importance in
estimating different invariants of the graph such as  girth, 
independence number and expansion coefficient. A large expansion coefficient
is determined by a small  $\lambda$ as shown in \cite{LPS1}. Lubotzky, Philips
and Sarnak, in \cite{LPS1}, have constructed a family of expander graphs called
Ramanujan graphs. Asymptotically, their graphs have the smallest possible
$\lambda$. They have called these expanders Ramanujan graphs, because all eigenvalues, except
the largest (of course in absolute value), satisfy Ramanujan's conjecture (or
more precisely Ramanujan-Petersson conjecture).\\ We shall give at first the
definition of regular Ramanujan hypergraphs. These hypergraphs are a natural generalization  of Ramanujan graphs. The eigenvalues of the adjacency
operators satisfy inequalities associated of a higher dimensional version of Ramanujan-Perterson conjecture.  In order to obtain a natural and simple
definition of Ramanujan hypergraphs, first we need some combinatorial
definitions and concepts. 
Our main references here are \cite{big} and \cite{bol}.
\begin{definition}
A graph is a pair of sets $(V_X,E_X)$ such that: \\ $E_X \subset \{\{x,y\}|
\{\{x,y\} \subset V_X\}$ and $ V_X \neq \varnothing$.\\
The set $V_X$ is the set of vertices of $X$ and $E_X$ is the set of edges of
$X$. The vertices x and y are said to be adjacent if $\{x,y\}$ is an edge. The
number of vertices adjacent to x is denoted by $d(x)$ and is said to be the
degree of x. If every vertex of X has degree s, then x is said to be
s-regular. \end{definition}
\begin{definition}
$i$If X is a graph with a finite number of vertices
$\{x_1,x_2,\ldots,x_n\}$, The adjacency matrix $A=[a_{ij}]$ of X is the
$n\times n$ matrix with entries $a_{ij}$ equal $1$ if $x_i$ is adjacent to $x_j$
and $0$ otherwise. \\
$ii)$ Denoting $L^2(V_X)$ $($or $L^2(X))$ the space of functions $f : V_X
\longrightarrow \mathbb{C}$ with the usual $L^2$-norm and with standard basis
the set of delta functions $\delta_v$, $v\in V_X$, the adjacency matrix
induces an operator \[ A : L^2(V_X) \longrightarrow L^2(V_X) \] with respect
to the basis $\{\delta | v \in V_X \}$.
 \end{definition}  
\begin{definition}
 A hypergraph $\bs{X}$ is a set $\bs{V}$ together with a family $\bs{\Sigma}$ of subsets of $\bs{V}$. The elements of $\bs{V}$ and $\bs{\Sigma}$  are called vertices and faces of the hypergraph. If $\bf{S} \in \bs{\Sigma}$, the rank of $\bf{S}$ is the cardinality $\bf{|S|}$ of $\bf{S}$ and the dimension of $\bf{S}$ is given by  $\bf{|S|}-1$.\end{definition}
We have seen the definition of labeled simplicial complex, and chamber
complex in chapter one. So we are ready to give the definition of labelable
hypergraph.
\begin{definition}
A chamber complex is a set $\Delta$, whose elements are called chambers,
together with a set ${\{\Large{\sim}}_i : i\in I\}$ of equivalence relations
on $\Delta$. Call $c,c'\in \Delta$ $i$-adjacent if $c {\Large{\sim}}_i c'$. In
this situation we refer to  $\Delta$ as  a chamber complex over
$I$\end{definition}
\begin{definition}
  Suppose that $c$ and $c'$ are in  $\Delta$, and that
  $c_0=c,c_1,\ldots,c_s=c'$ is a finite sequence of chambers such that
  $c_{k-1}$ is adjacent to and distinct from $c_k$ for each
  $k\in\{1,2,\ldots,s\}$. Then $(c_0,c_1,\ldots,c_s)$ is called a $\bs{gallery}$
  from $c$ to $c'$. If $c_{k-1}{\Large{\sim}}_{i_k}c_k$ for each
  $k\in\{1,2,\ldots,s\}$, then we say that the   $\bs{gallery}$ is of $\bs{type}$ $(i_1,i_2,\ldots,i_s)$
\end{definition} 
Very occasionally we need to consider slightly more general galleries, where
the requirement that $c_{k-1}\not =c_k$ for each $k$ is dropped. Such things will
 be called $\bs{stutter}$-$\bs{galleries}$.
\begin{definition}
 A $\bs{simplicial}$ $\bs{complex}$ is a set $X$ together with a set
 $\mathcal{S}$ of subsets of $X$ such that \begin{align*}
&1.& \text{each singleton subset}\quad \{x\}\quad \text{of}\quad
X\quad\text{belongs to}\quad \mathcal{S}.\\&2.& \text{if}\quad S\in
 \mathcal{S}\quad\text{and if}\quad A\subset S, then A\in  \mathcal{S};     
 \end{align*}
The sets $S\in \mathcal{S}$ are called $\bs{simplexes}$ or
$\bs{simplces}$. If $S\in \mathcal{S}$ and $A\subset S$, the simplex $A$ is
called a $\bs{face}$ of $S$. If $s$ is a simplex and $S$ is not a proper
subset of any simplex $S'$, we call $S$ a maximal simplex, or a $\bs{chamber}$.
\end{definition}
\begin{definition} \label{abcd123}
A simplicial complex $X$ is called  $\bs{labeled}$, if there is a set $I$ of
''labels'' or ''types'', so that each vertex $v$ has a type $\mathbb{t}(v)\in
I$, and such that each chamber has exactly one vertex of each type. In other
words, $\mathbb{t}:X\longrightarrow I$ is map such that for each chamber $C$,
the restriction $\mathbb{t}_{|_C}$ of $\mathbb{t}$ to $C$ is a bijection
$C\longrightarrow I$. This implies that each chamber has exactly $Card(I)$
vertices. \end{definition}
Under suitable conditions one can obtain labeled simplicial complex from a
chamber complex and vice versa. $($see \cite{Ro1} and   \cite{Ro2}$)$ for more
detail.\begin{remark} Specially our building $X_{\bs{\cdot}}(F^d)$ has satisfies all
   conditions in view of above correspondence between chamber and
  labeled simplicial complex.\end{remark}
\begin{definition}
A hypergraph X is labelable if it is a chamber  complex and there exist a set
I and a function which assigns to each vertex of X an element of I in such a
way that the vertices of every chamber are mapped bijectively onto
I. \end{definition}
\begin{definition}
Let $\bs{X}$ be a labelable hypergraph $\bs{X}$ with the label set
$\bs{I}=\{0,1,\dots,s\}$, if every vertex $x \in \bs{X}$ has exactly
$\bs{n_k}$  number of neighbor of type k, then $\bs{X}$ is called a
$\bs{(n_1,n_2, \dots,n_s)}$-regular hypergraph. In this situation we can
associate to any vertex $x$ of $\bs{X}$ a function $\mathbb{t}_x$ defined by 
$\mathbb{t}_x(y):= \mathbb{t}(x)-\mathbb{t}(y) \bmod d$. where $\ \mathbb{t}$
is the label map defined as the definition $(\ref{abcd123})$.
\end{definition}

denoted by $\underline{X}$.
\begin{definition}
Let $\bs{X}$ be a $\bs{(n_1,n_2,\dots,n_s)}$-regular hypergraph. We define for every $k \in \{0,1,\dots,s\}$  the $k^{th}$ adjacency matrix $\bs{A^{(k)}}$ as follows : denote
 for every $k \in \{0,1,\dots,s\}$ and for any two vertices x,y 
\begin{equation}\label{adja}
 \bs{\varepsilon^{(k)}}(x,y)=
                    \begin{cases}
                     1 & \qquad   \text{if $ \mathbb{t}_x(y)= k$}\\
                     0 & \qquad  \text{otherwise} 
                     \end{cases}
                    \end{equation} 
now $\bs{A^{(k)}}(i,j):=\bs{\varepsilon^{(k)}}(x_i,x_j)\bs{A}(i,j)$ where  $\bs{A}(i,j)$ is the $ij^{th}$ entries of the adjacency matrix of underlying graph $\underline{\bs{X}}$.
\end{definition}
 We are  ready now to give our main definition in this thesis, namely the
 definition of a Ramanujan hypergraph.\\
\begin{definition}\label{Ramali}
A $\bs{(n_1,n_2,\dots,n_{d-1})}$-regular hypergraph $\bs{X}$ is called a
Ramanujan hypergraph with the bound $\bs{(c_1,c_2, \dots,c_{d-1})}$ where for
all \\$k\in \{1,2,\ldots,d-1\}$, $\bs{c_k}$ are positive real number, if every eigenvalue $\bs{\lambda^{(k)}}$ of the  $k^{th}$ adjacency matrix $\bs{A^{(k)}}$ of $\bs{X}$ is either  $\bs{\lambda^{(k)}=n_k}$ or $|\bs{\lambda^{(k)}| \leqslant c_k}$.
\end{definition}
  As we will see in the next chapter the structure
$\Gamma\backslash vert(X_{\bs{\cdot}}(F^d)$ 
will be a finite regular hypergraph. These quotients would be our main object
of study in the next chapters.
\section{Quotient hypergraphs, the general situation}
 $F = \mathbb{F}_q(t)$ is again a rational function field over
 $\mathbb{F}_q$. $\{\underline{p}_1,\ldots, \underline{p}_m\}$ is a set of places of $F$, such that
 $p_m = \infty$, $\underline{p}$ is a prime of degree one of $F$, such that
 \[
 \underline{p} \notin \{\underline{p}_1,\ldots, \underline{p}_m\} \]
$p = p(t)$ is an irreducible polynomial, such that $\underline{p} =
(p(t))$. We have the orders $\mathbb{F}_q[t][\frac{1}{p}]$ as well as 
\[ \mathcal{O}_F^{(p)} : = \{ x \in  \mathbb{F}_q[t][\frac{1}{p}] \quad |
\quad v_\infty(x) \geqslant 0 \} \]
$D$ is  a central division algebra over $F$ of dimension $d^2$, which is
ramified exactly at  $\{\underline{p}_1,\ldots, \underline{p}_m\}$. For
simplicity we assume, that $D$ is totally ramified in $\infty$, such that
there exists a unique valuation $v_\infty$ of $D$ extending the corresponding
valuation of $F$. $\mathcal{O}_D$ is a maximal order of $D$ over
$\mathbb{F}_q[t]$. Then \begin{equation*}
\mathcal{O}_D [\frac{1}{p}] : =
\mathbb{F}_q[t][\frac{1}{p}]\otimes_{\mathbb{F}_q[t]}\mathcal{O}_D
\end{equation*}
is a maximal order in $D$ over $\mathbb{F}_q[t][\frac{1}{p}]$, similarly
\begin{equation*}
\mathcal{O}_D^{(p)} : =\{ x \in  \mathcal{O}_D[\frac{1}{p}] \quad |
\quad \tilde{v}_\infty(x) \geqslant 0 \}
\end{equation*}$($where $\widetilde{v}$ is standard extension of $v$ into $D)$
is a maximal order in $D$ over $\mathcal{O}_F^{(p)}$ .
\begin{definition} 
\begin{equation*}
   \Gamma(1) : = \biggl(\mathcal{O}_D
  [\frac{1}{p}]\biggr)^\times/ Z,
 \end{equation*}
where $Z$ is the center of the group ${\mathcal{O}_D[\frac{1}{p}]}^\times$, the
unit group of the maximal order  $\mathcal{O}_D[\frac{1}{p}]$
\end{definition}
We are considering in the following congruence subgroups $\Gamma$ of
$\Gamma(1)$. We have the embedding 
 \begin{equation*}
   \Gamma(1) \hookrightarrow D_{\underline{p}}^ \times/Z \cong \text{\rm{GL}}(d, \mathbb{F}_q(t)_p)/Z
 \end{equation*}
as obtained earlier and the corresponding action of the different groups on
the Bruhat-Tits building $X_{\bs{.}}({\mathbb{F}_q(t)_p}^d)$.\\
We consider now  only congruence subgroups $\Gamma \subset \Gamma(1)$, which
are torsion free. In particular their action on
$X_{\bs{.}}({\mathbb{F}_q(t)_p}^d)$ is fixpoint free. Again the quotient
complex $($or quotient hypergraph, if preferred$)$ is a locally labeled
complex resp. locally labeled hypergraph. We have quite generally :
\begin{theorem} \label{finite-hyper}The quotient complex 
\[ \Gamma \backslash  X_{\bs{.}}({\mathbb{F}_q(t)_p}^d) \] is a finite complex.
\end{theorem}
\begin{proof} This follows from Godement's criterion.
\end{proof}
\begin{remark} Of course here one can argue also directly, because the class
  number of a maximal order $\mathcal{O}_D$ $($and all the variants above$)$
  is finite. We consider the $\mathbb{C}$-vector space of $\mathbb{C}$-valued
  cochains: 
\begin{equation*}
C^0(\Gamma \backslash  X_{\bs{.}}({\mathbb{F}_q(t)_p}^d)) = : C^0(\Gamma
\backslash  X_{\bs{.}}) := \text{Map}(\Gamma
\backslash  X_0; \mathbb{C}).
\end{equation*}
 Because $ X_{\bs{.}}$ is a locally labeled simplicial complex, and the action
of $\Gamma$ is compatible with the local labeling, also $\Gamma \backslash
X_{\bs{.}}({\mathbb{F}_q(t)_p}^d)$ is a locally labeled simplicial
complex.\end{remark}
 \begin{definition}
For each $i, 1 \leqslant i \leqslant d-1$, we have the $i$-th adjacency
operator
\begin{align*}
A^{(i)} &: C^0(\Gamma \backslash  X_{\bs{.}}) \longrightarrow C^0(\Gamma
\backslash  X_{\bs{.}})\\
&f \longmapsto A^{(i)}(f)
\end{align*}
where \[ A^{(i)}(f)(x) = \sum_{t(y;x) = i}f(y) \]
\end{definition}
\begin{remark} We remind, that  $t(y;x) = i$ means:
\begin{equation*}
y \in \Gamma \backslash  X_0, \quad \text{such that} \quad <x,y> \ \in   \Gamma \backslash  X_1,
\end{equation*}
\end{remark}
that is, $x$ and $y$ are neighbors resp. there is a one-simplex joining them and
furthermore the type of the vertex $yA$ with respect to $x$ is $i$.\\
Next, we give an adelic description of the set $\Gamma \backslash  X_0$ above
and similarly of the associated space of $\mathbb{C}$-valued functions. Denote
$\mathbbb{A}_F$ the adeles over $F$, similarly
\begin{equation*}
 D^\times(\mathbb{A}_F) := \prod_{v\in |F|}(D^\times(F_v);
 D^\times(\mathcal{O}_v))
\end{equation*}
the adeles of the multiplicative group of $D$, $D^\times$. Z denotes in this
context again the center of $D^\times(\mathbb{A})$, we have $($ in the sense of
algebraic groups$)$ 
\begin{equation*}
 (D^\times/Z)(\mathbb{A}_F) = D^\times(\mathbb{A}_F)/Z
 \end{equation*}
The order $\mathcal{O}_D$ defines local orders $\mathcal{O}_{D,\underline{r}} \subset D_{\underline{r}}$
for all primes $r \neq \infty$. We have identified earlier
 \[ D_{\underline{p}} \quad
\widetilde{\longrightarrow} \quad \mathbb{M}(d, F_{\underline{p}}) \]
and thereby $\mathcal{O}^\times_{D,\underline{p}}/Z$ with the stabilizer of the standard
lattice 
\[L_0 = \mathcal{O}_{\underline{p}}^d = (\mathbb{F}_q[t]_p)^d \subset
F^d_{\underline{p}}. \] \begin{definition} \ \\ 
$i)$ We have \begin{equation*}
 \mathfrak{K}^{(1)} : = \prod_{\underline{r} \neq \infty,
   \underline{p}}(\mathcal{O}_{D,\underline{r}}^\times/Z) \times D_\infty^\times/Z ,
 \end{equation*}
 an open subgroup of $(D^\times(\mathbb{A}_F)/Z$.\\
$ii)$ Any congruence subgroup $\Gamma \subset \Gamma(1)$ of finite index
defines in a natural way a congruence subgroup  $\mathfrak{K}\subset \mathfrak{K}^{(1)}$
  of finite index.
 \end{definition}
\begin{proposition} \label{above1}There is a natural bijection 
 \begin{equation*}
 \Gamma \backslash X_0(\mathbb{F}_q(t)_p^d) \quad \widetilde{\longrightarrow}\quad
 D^\times \backslash D^\times(\mathbb{A}_F)/Z . \mathfrak{K} . D^\times(\mathcal{O}_{\underline{p}})
 \end{equation*}
\end{proposition}
\begin{proof} We can identify 
\begin{equation*}
 X_0(\mathbb{F}_q(t)_p^d) \quad \widetilde{\longrightarrow}\quad
  D^\times(F_{\underline{p}})/ D^\times(\mathcal{O}_{\underline{p}}).Z .
 \end{equation*}
Therefore we obtain a mapping 
\begin{equation*}
 X_0(\mathbb{F}_q(t)_p^d) \quad \longrightarrow \quad
 D^\times(\mathbb{A}_F)/Z . \mathfrak{K} . D^\times(\mathcal{O}_{\underline{p}})
\end{equation*}
This mapping induces a map 
\begin{equation*}
 \Gamma \backslash X_0(\mathbb{F}_q(t)_p^d) \quad \widetilde{\longrightarrow}\quad
 D^\times \backslash D^\times(\mathbb{A}_F)/Z . \mathfrak{K} . D^\times(\mathcal{O}_{\underline{p}}).
 \end{equation*}
We construct now an inverse map. Suppose, given an adele $(x_v) \in
D^\times(\mathbbb{A}_F)$, and thereby the doubled class 
 \begin{equation*}
D^\times(x_v)(Z.\mathfrak{K}.D^\times(\mathcal{O}_{\underline{p}})).
\end{equation*}
Using strong approximation theorem and also,that the class number of
$\mathbb{F}_q[t]$  is one, we can find $\gamma \in D^\times$, such that
\begin{equation*}
D^\times(\gamma x_v)_{v \in |F|}(Z.\mathfrak{K}.D^\times(\mathcal{O}_{\underline{p}})),
\end{equation*}
is of the form 
\begin{equation*}
D^\times(y_v)(Z.\mathfrak{K}.D^\times(\mathcal{O}_{\underline{p}})),
\end{equation*}
where \begin{equation*}
y_v = \begin{cases}
      1    & \text{for all}  \quad v \neq \underline{p}\\
      y_{\underline{p}} \in  D^\times(F_{\underline{p}}) &\text{at}\quad 
      \underline{p}
\end{cases}
 \end{equation*} 
We construct the inverse map then as 
\begin{equation*}
D^\times(\gamma x_v)_v (Z.\mathfrak{K}.D^\times(\mathcal{O}_{\underline{p}}))
\longmapsto \Gamma
. (y_{\underline{p}}).Z.D^\times(\mathcal{O}_{\underline{p}})) \in \Gamma
\backslash D^\times(F_{\underline{p}}) /
D^\times(\mathcal{O}_{\underline{p}}).Z .
\end{equation*}
It is immediate to see, that this gives a well defined map, where one has to
take into regard additionally, that we have, considering
$D^\times(F_{\underline{p}})$ as   a subgroup of $D^\times(\mathbbb{A}_F)$
\begin{equation}\label{cabali}
 D^\times(F_{\underline{p}}) \cap \mathfrak{K} = \Gamma
\end{equation}
and 
\begin{equation*}
 X_0(\mathbb{F}_q(t)_p^d) \quad \widetilde{\longrightarrow}\quad
  D^\times(F_{\underline{p}})/ D^\times(\mathcal{O}_{\underline{p}})Z .
 \end{equation*}
\end{proof}
We have therefore immediately
\begin{cor}
There is a canonical identification 
\begin{equation*}
C(\Gamma \backslash  X_0(\mathbb{F}_q(t)_p^d); \mathbb{C}) \quad
\widetilde{\longrightarrow}\quad C(D^\times \backslash
D^\times(\mathbb{A}_F)/Z . \mathfrak{K}
. D^\times(\mathcal{O}_{\underline{p}}); \mathbb{C}).
 \end{equation*}
In particular the Hecke-algebra of biinvariant functions
$\mathcal{H}(D^\times(\mathbb{A}_F)//  \mathfrak{K}
. D^\times(\mathcal{O}_{\underline{p}}))$ is acting on  $C(\Gamma \backslash
X_0(\mathbb{F}_q(t)_p^d); \mathbb{C})$, by convolution. As a subalgebra the
spherical Hecke algebra $\mathcal{H}(D^\times(F_{\underline{p}})//
 D^\times(\mathcal{O}_{\underline{p}}))$ is acting on $C(\Gamma \backslash
X_0(\mathbb{F}_q(t)_p^d); \mathbb{C})$.\end{cor}
\begin{proof}
The identification of the space of functions follows immediately from
Proposition \ref{above1}. The assertions concerning the action of the Hecke
algebra are also obvious, but there will be some discussion in $3.3.$
\end{proof}
\section{Adjacency operators and Hecke operators} 
In the last section we had identified the set of vertices 
\begin{equation*}
\Gamma \backslash  X_0(F_{\underline{p}}^d)
\end{equation*}
of the quotient hypergraph associated to the discrete group $\Gamma \subset
D^\times(F_{\underline{p}})/Z$ with the set of adelic double classes
\begin{equation*}
D^\times (F) \backslash D^\times(\mathbb{A}_F)/Z . \mathfrak{K}
. D^\times(\mathcal{O}_{\underline{p}}).
\end{equation*}
Correspondingly we have an identification of the spaces of functions 
\begin{equation*}
C(\Gamma \backslash X_0(\mathbb{F}_q(t)_p^d); \mathbb{C}) \quad \widetilde{\longrightarrow}\quad C( D^\times (F) \backslash D^\times(\mathbb{A}_F)/Z . \mathfrak{K}
. D^\times(\mathcal{O}_{\underline{p}}); \mathbb{C}).
\end{equation*}
The Hecke algebra associated to
$\mathfrak{K}.D^\times(\mathcal{O}_{\underline{p}})$ in
$D^\times(\mathbb{A}_F)$ is acting and in particular  the local spherical Hecke
algebra \[\mathcal{H}(D^\times(F_{\underline{p}});
Z.D^\times(\mathcal{O}_{\underline{p}})). \]
For the convenience of the reader we remind at this point about some basic
concepts of the theory of automorphic representations in connection with the
situation above. So, 
\begin{equation*}
C(D^\times(F)\backslash (\mathbb{A}_F)/Z . \mathfrak{K}. D^\times(\mathcal{O}_{\underline{p}}); \mathbb{C})
\end{equation*}
is a finite dimensional subspace of the space of functions
\begin{equation*}
C_c(D^\times(F)\backslash D^\times(\mathbb{A}_F)/Z ; \mathbb{C})
\end{equation*}
On this infinite dimensional $\mathbb{C}$-vector space there is an obvious
smooth action $R$ of the locally compact topological group
$D^\times(\mathbb{A}_F)$ by \[ R(x)(f)(g) = f(gx), \]
where $ g \in D^\times(F)\backslash D^\times(\mathbb{A}_F)/Z $,
$x \in  D^\times(\mathbb{A}_F) \quad \text{resp.}\quad
(D^\times/Z)(\mathbb{A}_F) \quad \text{(as wanted)}$ and $f \in
C_c^\infty(D^\times(F)\backslash D^\times(\mathbb{A}_F)/Z ; \mathbb{C})$
the space of locally constant functions on 
$D^\times(F)\backslash D^\times(\mathbb{A}_F)/Z $,  
with compact support. We have the identification of the subspace of
$\mathfrak{K}. D^\times(\mathcal{O}_{\underline{p}})$-invariants: 
\begin{equation*}
{C_c(D^\times(F)\backslash D^\times(\mathbb{A}_F)/Z ; \mathbb{C})}^{\mathfrak{K}. D^\times(\mathcal{O}_{\underline{p}})} \quad \widetilde{\longrightarrow} \quad C(D^\times(F)\backslash D^\times(\mathbb{A}_F)/Z . \mathfrak{K}. D^\times(\mathcal{O}_{\underline{p}}); \mathbb{C}).
\end{equation*}
There remains the action of the Hecke algebra of biinvariant functions
\begin{align*}
 \mathcal{H}(D^\times(\mathbb{A}_F)//  \mathfrak{K}
. D^\times(\mathcal{O}_{\underline{p}}))& \\ &
\cong (\otimes_{\underline{r}\neq
  \underline{p},\infty}\mathcal{H}(D^\times (F_{\underline{r}})//
  Z(\mathfrak{K}_{\underline{r}})) \otimes \mathcal{H}(D^\times(F_{\underline{p}})//
  D^\times(\mathcal{O}_{\underline{p}}))
\end{align*}
where 
\begin{equation*}
 \mathfrak{K} = (\prod_{\underline{r}\neq
  \underline{p},\infty}\mathfrak{K}_{\underline{r}} ) \times D^\times(F_{\infty})\times
  D^\times(\mathcal{O}_{\underline{p}}).
\end{equation*} The $($local$)$ spherical Hecke algebra 
\begin{equation*}
 \mathcal{H}(D^\times(F_{\underline{p}})// Z
.D^\times(\mathcal{O}_{\underline{p}})) \cong\\ 
\mathcal{H}(\text{\rm{GL}}(d;F_{\underline{p}})// Z
.\text{\rm{GL}}(d;\mathcal{O}_{\underline{p}}))
\end{equation*}
has as $\mathbb{C}$-vector space the generators  
\[
\text{\rm{GL}}(d;\mathcal{O}_{\underline{p}})\  \text{diag}(p^{n_1}, \ldots,
p^{n_d})\  
\text{\rm{GL}}(d;\mathcal{O}_{\underline{p}}).Z/Z
\end{equation*}
 where we have
\begin{equation*}
(n_1,\ldots,n_d) \in \mathbb{Z}^d/Z\bs{.} (1,\ldots, 1) ,\quad n_1 \leqslant
n_2 \leqslant \ldots \leqslant n_d.
\end{equation*}
 We denote the set
 \begin{equation*}
A_+ := \{(n_1,\ldots,n_d) \in \mathbb{Z}^d/Z\bs{.} (1,\ldots, 1)\quad |\quad n_1 \leqslant
n_2 \leqslant \ldots \leqslant n_d\}. 
\end{equation*}
\begin{proposition}
 As a $\mathbb{C}$-vector space
\begin{align*}
 &\mathcal{H}(D^\times(F_{\underline{p}})/ Z
.D^\times(\mathcal{O}_{\underline{p}})) \cong\\ 
&\bigoplus _{(n_1,\ldots,n_d) \in A_+}\biggl(\mathcal{H}(D^\times(\mathcal{O}_{\underline{p}}))
 \ \text{diag}(p^{n_1}, \ldots,
p^{n_d})\ \mathcal{H}(D^\times(\mathcal{O}_{\underline{p}})
)\bs{.} Z/Z\biggr)
\end{align*}
\end{proposition}
 The different Hecke algebras form, as the word indicates, algebras under the
convolution of functions, that is 
\begin{equation*}
(f_1\ast f_2)(x) := \int_{D^\times(F_{\underline{p}})/Z} f_1(xy^{-1} d)f_2(y) \mu(y)
\end{equation*}
 where $d\mu$ denotes the Haar measure on $D^\times(F_{\underline{p}})/Z$, In
particular, we have the  double classes
\begin{equation*}
D^\times(\mathcal{O}_{\underline{p}})
) \ \text{diag}(1,\ldots,1,p, \ldots,
p)\ D^\times(\mathcal{O}_{\underline{p}})
)\bs{.} Z/Z
\end{equation*}
 $(i$-times 1 above, at least one $p$ occurring$)$ and the associated
characteristic functions $\chi_{\underline{p},i}$.
\begin{definition}
 The Hecke operator $T_{\underline{p},i}$ is given as the convolution operator 
 \begin{align*}
  C(D^\times(F)\backslash D^\times(\mathbb{A}_F)/Z
 . \mathfrak{K}. D^\times(\mathcal{O}_{\underline{p}}); \mathbb{C}) \quad
 &\longrightarrow \quad C(D^\times(F)\backslash D^\times(\mathbb{A}_F)/Z
 . \mathfrak{K}. D^\times(\mathcal{O}_{\underline{p}}); \mathbb{C})\\
f \quad &\longmapsto f\ast \chi_{\underline{p},i} = : T_{\underline{p},i}(f).
\end{align*}  
\end{definition}
 We will now identify the Hecke operators $T_{\underline{p},i}$ with the
corresponding adjacency operators. The Hecke operators $T_{\underline{p},i}$
 can be seen best using the isomorphism 
 \begin{equation*}
 C(D^\times(F)\backslash D^\times(\mathbb{A}_F)/Z
 . \mathfrak{K}. D^\times(\mathcal{O}_{\underline{p}}); \mathbb{C}) \quad
 \longrightarrow \quad C(\Gamma \backslash D^\times(F_{\underline{p}})/Z. D^\times(\mathcal{O}_{\underline{p}}); \mathbb{C})
\end{equation*}
 As the Hecke operators $T_{\underline{p},i}$ are acting by convolution from
the right side, they will be commute with the action of $\Gamma$ from the left
 side. It is therefore enough to compute the action of the Hecke operators
 $T_{\underline{p},i}$ on the space of functions 
\begin{equation*}
C(D^\times(F_{\underline{p}})/Z. D^\times(\mathcal{O}_{\underline{p}});
\mathbb{C})\quad
\text{resp.}\quad
C^{\infty}(D^\times(F_{\underline{p}})/Z. D^\times(\mathcal{O}_{\underline{p}});
\mathbb{C}).
\end{equation*}
 As it will turn out, the Hecke operators will be defined locally and act
therefore on both spaces. We have the identification 
\begin{align*}
 D^\times(F_{\underline{p}})/Z. D^\times(\mathcal{O}_{\underline{p}})\quad
&\widetilde{\longrightarrow}\quad \text{\rm{GL}}(d; F_{\underline{p}})/ \text{\rm{GL}}(d;
    \mathcal{O}_{\underline{p}})\bs{.} Z\\
&\widetilde{\longrightarrow}\quad X_0(F_{\underline{p}}^d).
\end{align*}
 Therefore, the action of the Hecke algebra $\mathcal{H}(D^\times(F_{\underline{p}})// Z
\bs{.}D^\times(\mathcal{O}_{\underline{p}}))$ can be seen also on the space of
functions $C(X_0(F_{\underline{p}}^d); \mathbb{C})$ 
as well as  $C_c(X_0(F_{\underline{p}}^d); \mathbb{C}).$
\begin{theorem}
 $i)$ The Hecke operator $T_{\underline{p},i}$ equals the adjacency operator
 $A^{(i)}$ for each $ i =1, \ldots,d-1$. \\
 $ii)$  The Hecke algebra at $\underline{p}$ and the algebra generated by the
adjacency operators $A^{(i)}$ $(i =1, \ldots,d-1)$ act equally.  \\
 $iii)$ The action of both algebras commutes with the action of
 \[\text{\rm{GL}}(d;
 F_{\underline{p}})\cong  D^\times(F_{\underline{p}}).\]
\end{theorem}
\begin{proof}
$iii)$ is obvious. For the case of adjacency operators it follows directly
from the formulas for the $A^{(i)}$, acting on the space of functions $
C(X_0(F_{\underline{p}}^d); \mathbb{C})$. For the case of the Hecke algebra it
follows, because $\text{\rm{GL}}(d;
F_{\underline{p}})\cong  D^\times(F_{\underline{p}})$ is acting from the left
side, where as the Hecke operators are given by convolution with biinvariant
functions from the right.\\
 $i)$ As a $\text{\rm{GL}}(d;
F_{\underline{p}})\cong  D^\times(F_{\underline{p}})$-module,
$C_c(X_0(F_{\underline{p}}^d); \mathbb{C})$ is generated by the characteristic
function \[ \chi_{D^\times(\mathcal{O}_{\underline{p}})\bs{.}Z/Z} =\chi_{L_0},\]
under the identification above. But 

\begin{align*}
 \quad & \chi_{D^\times(\mathcal{O}_{\underline{p}})\bs{.}Z/Z} \ast T_{\underline{p},i}\\
 =\quad  & \chi_{D^\times(\mathcal{O}_{\underline{p}})\bs{.}Z/Z} \ast  \chi_{D^\times(\mathcal{O}_{\underline{p}})} \text{diag}(1,\ldots,1,p, \ldots,
p) \chi_{D^\times(\mathcal{O}_{\underline{p}})\bs{.} Z/Z}\\
=\quad & \chi_{D^\times(\mathcal{O}_{\underline{p}})} \text{diag}(1,\ldots,1,p, \ldots,
p)\ \chi_{D^\times(\mathcal{O}_{\underline{p}}) \bs{.}Z/Z}\\ 
= \quad & \sum_{\text{type}(L;L^{\prime}) = i}\chi_{L^{\prime}}\\
 = \quad  & A^{(i)}(\chi_L)
\end{align*}
Here we have used, that $\chi_{D^\times(\mathcal{O}_{\underline{p}})\bs{.}
  Z/Z}$ is the unit element of the Hecke algebra
  $\mathcal{H}(D^\times(F_{\underline{p}})// Z
  \bs{.}D^\times(\mathcal{O}_{\underline{p}}))$. This altogether shows
  $i)$. \\ $ii)$ is an immediate consequence of $i)$, as both algebras of
  operators are generated by the Hecke operators $T_{\underline{p},i}$  resp. the $A^{(i)}$ $(i =1, \ldots,d-1)$.
 \end{proof}
\section{The Ramanujan Property}
 As  is well known, one can decompose the representation space 
\begin{equation*}
C_c(D^\times(F) \backslash D^\times(\mathbbb{A}_F) / Z; \mathbb{C}) \quad = \quad
\bigoplus_{\pi} V_{\pi}
\end{equation*}
into a direct sum of irreducible automorphic representations, $(V_{\pi}, \pi)$
of $D^\times(\mathbbb{A}_F)$. The subspace 
\begin{equation*}
C(D^\times(F) \backslash D^\times(\mathbbb{A}_F) /
Z.\mathfrak{K}. D^\times(\mathcal{O}_{\underline{p}}); \mathbb{C}) \subset C(D^\times(F) \backslash D^\times(\mathbbb{A}_F) / Z; \mathbb{C})
\end{equation*} can be written then in the form 
\begin{equation*}
   C(D^\times(F) \backslash D^\times(\mathbbb{A}_F) /
Z.\mathfrak{K}. D^\times(\mathcal{O}_{\underline{p}}); \mathbb{C})\end{equation*}\begin{equation*}
= \bigoplus_{\pi}(V_{\pi})^{\mathcal{H}(D^\times(\mathbbb{A}_F)//Z.\mathfrak{K}.D^\times(\mathcal{O}_{\underline{p}}))}
\end{equation*}
\begin{equation*}
= \bigoplus_{\pi}(V_{\pi})^{\mathcal{H}(D^\times(\mathbbb{A}_F^{(p)})//Z.\mathfrak{K})}\otimes(V_{\underline{p}}^{D^\times(\mathcal{O}_{\underline{p}})})
\end{equation*}
decomposing the representation \[ V_{\pi} = ( V_{\pi^{\underline{p}}} \otimes
V_{\underline{p}}) \]
and correspondingly 
\begin{equation*}
(V_{\pi})^{\mathcal{H}(D^\times(\mathbbb{A}_F)//Z.\mathfrak{K}.D^\times(\mathcal{O}_{\underline{p}}))}
\end{equation*}
\begin{equation*}
\cong \quad
(V_{\pi^{(\underline{p})}}^{\mathcal{H}(D^\times(\mathbbb{A}_F^{(p)})//Z.\mathfrak{K})})\otimes
  (V_{\underline{p}}^{D^\times(\mathcal{O}_{\underline{p}})}).
\end{equation*}
Therefore, only those representations $\pi$ occur in the decomposition, for
which \[   V_{\underline{p}}^{D^\times(\mathcal{O}_{\underline{p}})} \neq 0\]
If this is the case, $V_{\underline{p}}$ is a spherical representation. As the
  local Hecke algebra $\mathcal{H}(D^\times(F_{\underline{p}})
  //D^\times(\mathcal{O}_{\underline{p}}))$ is commutative, it follows that
  $V_{\underline{p}}^{D^\times(\mathcal{O}_{\underline{p}})}$ has to be
  one-dimensional. The action of \[\mathcal{H}(D^\times(F_{\underline{p}})
  //D^\times(\mathcal{O}_{\underline{p}})),\] is then given by a character $($
  ring homomorphism $)$ \[ \chi_{\pi_{\underline{p}}} :
  \mathcal{H}(D^\times(F_{\underline{p}})//D^\times(\mathcal{O}_{\underline{p}}))
  \quad \longrightarrow \quad \mathbb{C} \] into the field of complex
  numbers. In particular, \[  \chi_{\pi_{\underline{p}}}(
  T_{\underline{p},i})\in \mathbb{C} \]
 and the Ramanujan property of the automorphic representation at
 $\underline{p}$  deals with these. We have to make one further remark:\\
Obviously in \[ C^{\infty}(D^\times(F) \backslash D^\times(\mathbbb{A}_F) /Z;
\mathbb{C}) ,\]
we have the subrepresentations given by characters $($multiplicative
homomorphisms$)$ \[ \chi : D^\times(F) \backslash D^\times(\mathbbb{A}_F) /Z
\longrightarrow \mathbb{C}^\times. \]
In particular the trivial character, corresponding to the subspace of constant
functions, occurs as $D^\times(\mathbbb{A}_F)$-invariant subspace.
\begin{definition} \label{Ramali}  $i)$ An irreducible automorphic representation $\pi$,
  such that  $\pi_{\underline{p}}$ is a spherical representation, is said to
  satisfy the Ramanujan property at $\underline{p}$, iff the eigenvalues  
 $\chi_{\pi_{\underline{p}}}( T_{\underline{p},i})$ satisfy
\[ |\chi_{\pi_{\underline{p}}}( T_{\underline{p},i})| \leqslant
  \binom{d}{i}q^{\frac{i(d-i)}{2}} |\sigma_i(z_1,\ldots,z_d)|,\]
where $z_1,\ldots,z_d \in \mathbb{C}$ are complex numbers with absolute value
$|z_j| = 1$ for $j = 1,\ldots, d$ and $\sigma_i(z_1,\ldots,z_d)$ $i$-th
elementary symmetric polynomial.\\
 $ii)$ The Ramanujan property at the prime $\underline{p}$ holds for the
 representation \[ C^{\infty}(D^\times(F) \backslash D^\times(\mathbbb{A}_F) /Z;
\mathbb{C}), \] if it holds in the sense above for all automorphic
 representations $\pi$ complementary to the invariant subspaces, generated by
 the multiplicative characters.\end{definition}
\begin{theorem} \label{main1}(W.Li.) The Ramanujan property at  $\underline{p}$ holds in
 the sense of Definition \ref{Ramali} $(i)$.\end{theorem}
\begin{cor} \label{cor47}The quotient hypergraphs $\Gamma \backslash X_0(F_p^d)$ are
  Ramanujan hypergraphs $($simplicial complexes$)$ in the sense of Definition
  \ref{Ramali} in Chapter $2$. 
\end{cor}
\begin{remark} In the case $d = 2$, the Ramanujan property can be shown by
  using the so called Jacquet- Langlands correspondence between automorphic
  representations of $D^\times(\mathbbb{A}_F)$ $($then $D$ is a quaternion
  algebra$)$ and automorphic representations of $\text{\rm{GL}}(2;\mathbbb{A}_F)$,
  not given by multiplicative characters. For these, the Ramanujan
  property is a consequence of results of Drinfeld \cite {Drin1}. Though in our situation
  one has now the recent results of L.Lafforgue, showing the Langlands
  correspondence for $\text{\rm{GL}}(d;\mathbbb{A}_F)$ for arbitrary $d$ and in
  particular the Ramanujan property for cuspidal automorphic representations,
  one can not conclude immediately here the Ramanujan property. What is missing, is
  a completely worked out Jacquet- Langland correspondence as above for $d =
  2$. Nevertheless, due to a trick from L.Clozel, one can conclude the
  Ramanujan property by working with the moduli scheme of $\mathcal{D}$-elliptic modules \cite{LRS},
  instead of working
  with the moduli scheme of shtukas. Concerning this, then reader has to
  consult \cite{wlij}.   We close this chapter indicating another method to
  show the Ramanujan property. We use a theorem of Arthur and Clozel
  \cite [Theorem $4.2.$]{Artur}, which gives the Jacquet-Langlands
  correspondence for the case $d$ a prime number. After that we can use the
  recent result of Lafforgue to conclude the Ramanujan property for this
  case. It might be mentioned that we found this approch at a time, when the
  preprints  \cite{wlij} and \cite{LSV1},  \cite{LSV2} avilable to us.
\end{remark}
Assume d is a prime number, $f$ is an irreducible element of $\mathbb{F}_q[t]$
with $(\underline{f}) \neq (\underline{p})$ and different from all other
primes in which $ D= \mathbb{F}_{q^d}(\tau)$ is ramified . The canonical homomorphism
\[\mathbb{F}_q[t]\frac{1}{p}] \longrightarrow \mathcal{O}_D[\frac{1}{p}]/f\mathbb{F}_q[t] [\frac{1}{p}] ,\] 
induces following homomorphism :\[\mathcal{O}_D[\frac{1}{p}] \longrightarrow \mathcal{O}_D[\frac{1}{p}]/f\mathcal{O}_D[\frac{1}{p}] ,\]
where as before $\mathcal{O}_D=
\mathbb{F}_{q^d}\{\tau\}$. So we obtain group homomorphism \[ \alpha_f^{(p)}: \Gamma(1)
\longrightarrow
(\mathcal{O}_D[\frac{1}{p}]/f\mathcal{O}_D[\frac{1}{p}])^\times /Z.\]
Recall that
$ \Gamma(1) = (\mathcal{O}_D[\frac{1}{p}])^\times /Z$. We define $\Gamma_f^{(p)} :=
\ker{\alpha_f^{(p)}}$. So  $\Gamma_f^{(p)}$ is a normal subgroup of $
\Gamma(1)$ of finite index and we have:
\begin{theorem}\label{abstract} $($Main Abstract Theorem$)$\\
$(1)$ $\Gamma_f^{(p)} \backslash X_{\bs{.}}(\mathbb{F}_q(t)_p^d)$ is a finite
$(n_1,\ldots,n_{d-1})$-regular graph with 
\begin{equation*}
n_i = \text {number of  $i$-dimensional
subspaces of} \  \mathbb{F}_{\widetilde{q}} \ \text{where}\ \widetilde{q} = q^{\deg{p}},
\end{equation*}for $i= 1,\ldots, d-1$, more precisely \[ n_i =
\binom{d}{i}_{\widetilde{q}} := \frac{\prod_{m = d- i +1}^d( {\widetilde{q}}^m -
  1)}{\prod_{m = 1}^i( {\widetilde{q}}^m - 1)} \ .\]
$(2)$ $\Gamma_f^{(p)} \backslash X_{\bs{.}}(\mathbb{F}_q(t)_p^d)$  is a
Ramanujan hypergraph in the sense of Definition \ref{Ramali} from Chapter
 $2$,i.e. It is Ramanujan with the bound $(c_1,\ldots,c_{d-1})$ where $c_i =
 \binom{d}{i} q^{\frac{i(d-i)}{2}\deg{p}}$ for $i = 1,\ldots, d-1$.
\end{theorem}
\begin{proof}
By Theorem \ref{finite-hyper} the quotient complex $($hypergraph$)$
$\Gamma_f^{(p)} \backslash X_{\bs{.}}(\mathbb{F}_q(t)_p^d)$ is a finite. the
expression about regularity is inherited from the structure of the Bruhat-Tits
Building $X_{\bs{.}}(\mathbb{F}_q(t)_p^d)$. So $(1)$ is done.
Following isomorphisms are  known from Chapter $3$:
\[ C(D^\times \backslash
D^\times(\mathbbb{A}_F)/Z.\mathfrak{K}.D^\times(\mathcal{O}_{\underline{p}}); \mathbb{C})
\ \widetilde{\longrightarrow} \ C(\Gamma 
\backslash D^\times(F_{\underline{p}})/ Z.D^\times(\mathcal{O}_{\underline{p}});
\mathbb{C})
\]
\[ C(\text{\rm{GL}}(d; F_{\underline{p}})/\text{\rm{GL}}(d; \mathbb{F}_q[t]_p).Z; \mathbb{C})
\ \widetilde{\longrightarrow} \ C(\text{PGL}(d; F_{\underline{p}})/\text{PGL}(d;  \mathbb{F}_q[t]_p); \mathbb{C})
\]
\[ \widetilde{\longrightarrow} \quad C(X_0(F_{\underline{p}}^d);  \mathbb{C}). \]
which hold for all congruence subgroups of $\Gamma(1)$. in particular for
$\Gamma_f^{(p)}$. Thus we define first associated to $f$ the congruence
subgroup \[ J_f :=\ker{(\mathcal{O}_{D,f}^\times/Z \longrightarrow
  (\mathcal{O}_{D,f}/f\mathcal{O}_{D,f})^\times /Z)} ,\]
and we define 
\[ \mathfrak{K} := \prod_{r\neq p,\infty,f} (\mathcal{O}_{D,\underline{r}}^\times
/Z). D^\times /Z(F_{\underline{p}}).D^\times /Z(F_{\infty}). J_f ,\]
and let $\mathfrak{M} := D^\times /Z. \mathfrak{K}$
\begin{remark} For any group $G$ here, we use notation $G^{(1)}$ for the
  subgroup of $G$ of elements with reduced norm $1$.\end{remark}
Applying Strong approximation Theorem two times sequentially, we see that
$D^\times/Z(\mathbbb{A}). F_{\mathbbb{A}}$ $( F_{\mathbbb{A}}$, diagonal
embedding of $F$ in $\mathbbb{A})$ is a finite index subgroup of
$\mathfrak{M}$. This plus \[ |  D^\times/Z(\mathbbb{A}) /
D^\times/Z(\mathbbb{A}). F_{\mathbbb{A}} | < \infty \ ,\]
shows $ |  D^\times /Z(\mathbbb{A}). F_{\mathbbb{A}}
/ \mathfrak{M} | < \infty $.  Also \[ [D^\times /Z \backslash
D^\times /Z(\mathbbb{A})/\mathfrak{K} :   D^\times /Z \backslash \mathfrak{M} /
\mathfrak{K}] < \infty \ . \]
But as we have seen by \ref{cabali}
\begin{equation*}
  \Gamma_f^{(p)} = D^\times(F_{\underline{p}}) \cap \mathfrak{K} .
\end{equation*}
So $ D^\times/Z \backslash \mathfrak{M} /
\mathfrak{K} \cong \Gamma_f^{(p)}\backslash D_{\infty}^\times /Z
D^\times(F_{\underline{p}})/Z$ and finally 
  \[ [D^\times/Z \backslash
D^\times/Z(\mathbbb{A})/\mathfrak{K} : \Gamma_f^{(p)}\backslash
D_{\infty}^\times  /Z
.D^\times (F_{\underline{p}})/Z  ] < \infty \ . \]
Let $\pi_\infty \otimes \pi_{\underline{p}}$ be an irreducible representation
of the right regular representation of $D^\times /Z(F_{\underline{p}})/
Z.D^\times /Z(F_{\underline{p}})$ in $C(\Gamma_f^{(p)}\backslash D^\times /Z(F_{\underline{p}})/
Z.D^\times /Z(F_{\underline{p}})$. There is an irreducible subrepresentation
of $D^\times /Z (\mathbbb{A})_F)$ in 
$\widetilde{\pi} = \otimes_r \widetilde{\pi}_{\underline{p}}$ such that
  $\widetilde{\pi}_{\underline{p}} = \pi_{\underline{p}}$  and
  $\widetilde{\pi}_{\infty} = \pi_{\infty}$.  \\
Suppose $\pi_{\underline{p}}$ does not occur in $V_{\underline{p}}^{D^\times /Z(\mathcal{O}_{\underline{p}}}
$ $($,i.e. It is not one-dimensional$)$. So $\widetilde{\pi}$ occurs in
$C(D^\times /Z \backslash D^\times /Z(\mathbb{A}_F); \mathbb{C})$. Since d is
assumed a prime number, we can apply the Weak Lifting Theorem \cite [Theorem
$4.2.$]{Artur}, and obtain a cuspidal subrepresentation  $\rho_\infty \otimes
\rho_{\underline{p}}$ of $ C(\text{PGL}(d; F_{\underline{p}})/\text{PGL}(d;
\mathbb{F}_q[t]_p); \mathbb{C})$, which is cuspidal and $\rho{\underline{p}} =
\pi{\underline{p}}$. By recent reslt of L. Lafforgue \cite{Laf3} we have \[ |
T_{\underline{p},i} | < \binom{d}{i} q^{\frac{i(d-i)}{2}\deg{p}}
\sigma_i(z_1',\ldots,z_d') \ ,\] where $z_1',\ldots,z_d' \in \mathbb{C}$ are complex numbers with absolute value
$|z_j'| = |\frac{z_j}{q^{(d-1)/2}}| = 1$ for $j = 1,\ldots, d$ and $\sigma_i(z_1',\ldots,z_d')$ $i$-th
elementary symmetric polynomial. Choose suitable $z_i$ as in Lafforgue's
expression of Ramanujan-Peterson conjecture ,i.e. the proof of $2$ is complete. 
\end{proof}
\section{The  Skew polynomial ring $\mathbb{F}_{q^d}\{\tau\}$}
In this section we collect various well known fundamental facts concerning skew
polynomial rings. These rings are well known in the theory of non commutative
rings and many of the facts we note below hold in greater
generality. However we have written up the relevant properties in the form
 we will need later on. For more details see 
 \cite{Carlitz1},\cite{Carlitz2},\cite{Carlitz3},\cite{Carlitz4},\cite{Carlitz5},\cite{Carlitz6} and  \cite{Ore1}, \cite{Ore2}. \\
We consider the finite field $\mathbb{F}_q$ of $q = l^n$ elements of
characteristic $l$. $\mathbb{F}_{q^d}$ is a finite extension of $\mathbb{F}_q$
of degree $d$.\\
We will construct now the skew polynomial ring  $\mathbb{F}_{q^d}\{\tau\}$. As a
set this is given by \begin{equation*}
\{\sum_{i=0}^n a_i \tau^i \quad | n\geq 0,\quad a_i \in \mathbb{F}_{q^d} \}
\end{equation*}
Addition is defined by \begin{equation*}
\sum_{i=0}^n a_i \tau^i \quad + \quad \sum_{i=0}^n b_i \tau^i :=
\sum_{i=0}^n (a_i + b_i) \tau^i .
\end{equation*}
For $\mathbb{F}_{q^d}\{\tau\}$ one obtains the structure of an
infinite dimensional vector space over $\mathbb{F}_q$.\\
Regarding the multiplication, the fundamental rule is\[ \tau \lambda =
\lambda^q \tau\]
for $\lambda \in \mathbb{F}_{q^d}$.  There is  a
unique multiplication  on  the ring $\mathbb{F}_{q^d}\{\tau\}$ satisfying this rule. 
\begin{proposition}
 The center of the ring   $\mathbb{F}_{q^d}\{\tau\}$ is given as
$\mathbb{F}_q[\tau^d]$.\end{proposition}
\begin{proof}
 It is obvious, that $\mathbb{F}_q[\tau^d]$ is contained in the center
$Z(\mathbb{F}_{q^d}\{\tau\})$. On the other hand, if $c = \sum_{i=0}^n a_i
\tau^i $ is a central element in  $\mathbb{F}_{q^d}\{\tau\}$, it has to commute
with all elements $\lambda \in \mathbb{F}_{q^d}$. This forces the $a_i$ for
$i \not \equiv 0 (\mod{d})$ to be zero. Therefore $c$ is of the form 
\begin{equation*} 
c= \sum_{\begin{subarray}{l}
          i \ \equiv \ 0 \mod{(d)}\\ 0 \leq i \leq n
            \end{subarray}}
             a_i(\tau^d))^{d/i}
\end{equation*}
Because  $\tau c = c \tau$, it follows additionally, that the $a_i$ occurring are
elements in  $\mathbb{F}_q$.\end{proof}
\begin{remark} It is immediate, that
 \begin{equation*}
 \{\alpha^i\tau^j \quad | \quad 0 \leq i,j \leq d-1\}
\end{equation*}
where \[\mathbb{F}_{q^d} = \mathbb{F}_q(\alpha)  \]
is a basis of the left (-right) modules $\mathbb{F}_{q^d}\{\tau\}$  over $\mathbb{F}_q[\tau^d]$.
\end{remark}
We denote $\tau^d =: t$. The center $\mathbb{F}_q[t]$ is of
  course the polynomial ring in the indeterminate $t$ over the finite field
  $\mathbb{F}_q$.
The following proposition is due to Ore, see \cite{Ore1}, \cite{Ore2}.
\begin{proposition} The skew polynomial ring $\mathbb{F}_q^d\{\tau\}$ is a
  left resp. right Euclidean ring.
\begin{proof}
One has to show the following property: given polynomials $f(\tau), g(\tau)\in
\mathbb{F}_{q^d}\{\tau\}$, there are polynomials $s(\tau)$ and  $r(\tau)$ such
that \begin{equation*}
f(\tau) =  s(\tau)g(\tau) +  r(\tau) \quad\text{such that}
\deg_{\tau}(r(\tau)) < \deg_{\tau}(g(\tau))
\end{equation*}
for the $\tau$-degrees of the  polynomials above in the obvious sense.
 This would show that $\mathbb{F}_{q^d}\{\tau\}$ is a left Euclidean ring.  But
 of course this can be seen by the usual division procedure of polynomials
 taking only into consideration, 
 that one can write $\lambda \tau = \tau \lambda^{1/q}$ because $\mathbb{F}_{q^d}
 \longrightarrow \mathbb{F}_{q^d}$,  $u \longmapsto u^q$ is a bijection. The
 property, that the ring  $\mathbb{F}_{q^d}\{\tau\}$ is right Euclidean means,
 one can find in the situation above elements  $s^\prime(\tau)$ and  $r^\prime(\tau)$ such
that \begin{equation*}
f(\tau) =  g(\tau)s^\prime(\tau) +  r^\prime(\tau) \quad\text{where}\\
\deg_{\tau}(r^\prime(\tau)) < deg_{\tau}(g(\tau))\quad\text{holds}.
\end{equation*}
This is shown in the same way.
\end{proof}
\end{proposition}
\begin{remark}
 All of this can be found in the literature $($as given above$)$
  even for
  more general skew polynomial rings $k\{\tau\}$, where $k$ is commutative field
  with automorphism $ \sigma : k \longrightarrow k$, such that the rule \[
  \tau \lambda = \sigma(\lambda) \tau\] holds for $\lambda \in k$.
\end{remark}
\begin{cor}
$i)$ Any left ideal $I$ in $\mathbb{F}_{q^d}\{\tau\}$ is a principal ideal of the
form $ I = \mathbb{F}_{q^d}\{\tau\}. a$, with $a\in I$ appropriate.\\
$ii)$ Similarly any right ideal $J$ in $\mathbb{F}_{q^d}\{\tau\}$ is a principal
ideal of the form  $ J = b.\mathbb{F}_{q^d}\{\tau\}$ with  $b \in J$ appropriate.
\begin{proof}
Given $I$, choose $g(\tau) \in I$ of minimal degree, if $I \neq
(0)$. Otherwise we are ready.  Suppose, $f(\tau)$ is an arbitrary element of
$I$. We can find  $s(\tau)$ and  $r(\tau)$ such
that \begin{equation*}
f(\tau) =  s(\tau)g(\tau) +  r(\tau)\\
 \quad \text{where}
\deg_{\tau}(r(\tau)) < deg_{\tau}(g(\tau)).
\end{equation*}
Because $f(\tau), g(\tau) \in I$, it follows that $r(\tau) \in I$, because $I$
is a left ideal. But then $r(\tau) = 0$, because otherwise we would have a
contradiction to the choice of $g(\tau)$ as nonzero element of $I$ with $
deg_{\tau}(g(\tau))$ minimal. Therefore we have  $f(\tau) = s(\tau)
g(\tau)$.\\
This implies that   $I =  \mathbb{F}_{q^d}\{\tau\}g(\tau)$ and therefore $I$ is
a principal left ideal.\\
$ii)$ is shown in the same way.
\end{proof}
\end{cor}
\begin{remark}
The generating elements as $g(\tau)$ above in i) are uniquely determined up to
multiplication by an element of $\mathbb{F}_{q^d}^\times$ from the left side.
\end{remark}
We have the following well known structure theorem for finitely generated
$\mathbb{F}_{q^d}\{\tau\}$-modules (left or right modules). It corresponds to
similar results for finitely generated modules over commutative principal
ideal domains:
\begin{theorem}
$i)$ Any finitely generated left $\mathbb{F}_{q^d}\{\tau\}$-module $M$ is the
direct sum of cyclic left $\mathbb{F}_{q^d}\{\tau\}$-modules:
\begin{equation*}
 M = \bigoplus_{i=1}^r \mathbb{F}_{q^d}\{\tau\}/\mathbb{F}_{q^d}\{\tau\} f_i(\tau)
\end{equation*}
where we can assume additionally 
\begin{equation*}
f_1(\tau) |_lf_2(\tau)|_l\ldots |_lf_r(\tau).
\end{equation*}
Here $a |_l b$ means left-divisibility, i.e. there is $r \in
\mathbb{F}_{q^d}\{\tau\}$ satisfying $r.a = b$.\\
$ii)$ If $M$ is a finitely generated torsion free left
$\mathbb{F}_{q^d}\{\tau\}$-module, then $M$ is a free
$\mathbb{F}_{q^d}\{\tau\}$-module, that is, the $f_i(\tau)$ above are all zero.
\end{theorem} 
\begin{proof}
see \cite{Carlitz1} or take any proof for the corresponding statement in  the
commutative situation.
\end{proof}
\begin{proposition}
$i)$ The ring $\mathbb{F}_{q^d}(\tau): =
\mathbb{F}_{q^d}\{\tau\}\bigotimes_{\mathbb{F}_q[t]}\mathbb{F}_q(t)$ ,
obtained by extension of the center $\mathbb{F}_q[t]$ of $\mathbb{F}_{q^d}\{\tau\}$, is as a left$($right$)$
module free of rank $d^2$ over the rational field $\mathbb{F}_q(t)$.\\
$ii)$ The center of $\mathbb{F}_{q^d}(\tau)$ is $\mathbb{F}_q(t)$ under the canonical
embedding of $\mathbb{F}_q(t)$ into $\mathbb{F}_{q^d}(\tau)$.\\
$iii)$ $\mathbb{F}_{q^d}(\tau)$ is a division algebra.  In
particular $\mathbb{F}_{q^d}(\tau)$ is a central simple algebra of dimension
$d^2$ over the rational function field   $\mathbb{F}_q(t)$.
\end{proposition}
\begin{proof}
$i)$ is clear, because $\mathbb{F}_{q^d}\{\tau\}$ is a free module of rank $d^2$
over  $\mathbb{F}_q[t]$.\\
$ii)$  If $u \in Z(\mathbb{F}_{q^d}(\tau))$, then there exists a nonzero
polynomial $f(t) \in \mathbb{F}_q[t]$ such that $f(t).u \in
Z(\mathbb{F}_{q^d}\{\tau\}) =  \mathbb{F}_q[t]$.\\
Conversely,  $\mathbb{F}_q(t)$, canonically embedded, is in the center of 
 $\mathbb{F}_{q^d}(\tau)$,\\
$iii)$ To show, that $\mathbb{F}_{q^d}(\tau)$ is a division algebra, consider the
canonical homomorphism of algebras 
\begin{align*}
\mathbb{F}_{q^d}(\tau)\bigotimes_{\mathbb{F}_q(t)}\mathbb{F}_q((t))
\longrightarrow \mathbb{F}_{q^d}((\tau))\\
\text{with} \quad f \otimes g  \longmapsto \quad f.g
\end{align*}
Here we are using  following notations:
\begin{equation*}
\mathbb{F}_q((t)) = \{ \sum_{i=N}^\infty a_i t^i\quad | \quad N\in \mathbb{Z}
\quad a_i \in \mathbb{F}_q\}
\end{equation*}
is the field of Laurent series $($at the place $t = 0$ of the field
$\mathbb{F}_q(t) )$,  $\mathbb{F}_{q^d}((\tau))$  is given as the skew field of
Laurent series 
\begin{equation*}
\mathbb{F}_q((\tau)) = \{ \sum_{i=N}^\infty a_i \tau^i\quad | \quad N\in \mathbb{Z}
\quad a_i \in \mathbb{F}_{q^d}\}
\end{equation*}

where we have again the communication rule \[ \tau \lambda = \lambda^q\tau\]
for $\lambda \in  \mathbb{F}_{q^d}$. It is immediate to see that
$\mathbb{F}_{q^d}((\tau))$ is a skew field. It is a $d$-dimensional vector
space over the field $\mathbb{F}_q((t))$ of Laurent series over $\mathbb{F}_q$.\\
Now, as vector space over $\mathbb{F}_q((t))$ both
$\mathbb{F}_{q^d}(\tau)\bigotimes_{\mathbb{F}_q(t)}\mathbb{F}_q((t))$ and
$\mathbb{F}_{q^d}((\tau)) $ are $d^2$-dimensional. Furthermore, the canonical
homomorphism 
\begin{equation*}
\phi : \mathbb{F}_{q^d}(\tau)\bigotimes_{\mathbb{F}_q(t)}\mathbb{F}_q((t))
\longrightarrow \mathbb{F}_{q^d}((\tau))
\end{equation*}
is surjective, as the image contains the elements $\{\alpha^i\tau^j \quad |
\quad 0 \leq i,j \leq d-1\}$, where $\mathbb{F}_{q^d} =\mathbb{F}_q(\alpha)$ as
above, which form a basis of the vector space  $\mathbb{F}_q^d((\tau))$ over
$\mathbb{F}_{q^d}((t))$. Then, as a surjective homomorphism between vector
spaces of equal dimension, $\phi$ is an isomorphism. Furthermore, $\phi$ is
compatible with the multiplicative structure, therefore $\phi$ is even an
isomorphism of algebras. Because $\mathbb{F}_{q^d}((\tau))$ is a skew field, the
algebra $\mathbb{F}_{q^d}(\tau)\bigotimes_{\mathbb{F}_q(t)}\mathbb{F}_q((t))$ is
a skew field as well. So the algebra $\mathbb{F}_{q^d}(\tau)$ $($as a sub algebra$)$
has no zero divisors. Furthermore it is a finite dimensional algebra over its
center $\mathbb{F}_q(t)$, which is a field. Then $\mathbb{F}_{q^d}(\tau)$ is a
skew field itself. This shows $iii)$
\end{proof}
\begin{remark}
 $i)$It follows from the considerations above, that the division algebra
 $\mathbb{F}_{q^d}(\tau)$ over its center $\mathbb{F}_q(t)$, where $t = \tau^d$,
 is ramified at the place $t=0$ with completion the skew field of Laurent
 series $\mathbb{F}_{q^d}((\tau))$,\\
$ii)$ Similarly at the place $t = \infty$, the place corresponding to the
 degree valuation of  $\mathbb{F}_q(t)$, $\mathbb{F}_{q^d}(\tau)$ is totally ramified.
\end{remark}
 We will show now, that these places are the only places of
 $\mathbb{F}_{q^d}(\tau)$ over  $\mathbb{F}_q(t)$, which are ramified.\\
The skew field  $\mathbb{F}_{q^d}(\tau)$ can be described as a cyclic algebra
 over its center  $\mathbb{F}_q(t)$ in the following way. First we have the
 cyclic Galois extension $\mathbb{F}_{q^d}(t)$ of $\mathbb{F}_q(t)$ with
 canonical embedding $\mathbb{F}_{q^d}(t) \hookrightarrow \mathbb{F}_{q^d}(\tau)$
 by mapping $t \mapsto t=\tau^d$. Of  course
 \begin{equation*}
Gal(\mathbb{F}_{q^d}(t) / \mathbb{F}_q(t) ) \cong Gal(\mathbb{F}_{q^d} /
\mathbb{F}_q) \cong <\overline{\tau}>
 \end{equation*}
where $\overline{\tau} :  \mathbb{F}_{q^d} \longrightarrow \mathbb{F}_{q^d},\quad
u \mapsto u^q$ is given by the Frobenius automorphism. $\mathbb{F}_{q^d}(t)$
is an unramified field extension of $\mathbb{F}_q(t)$ at all  places of
$\mathbb{F}_q(t)$. As it can be  compute the invariants of the skew
field   $\mathbb{F}_{q^d}(\tau)$ over $\mathbb{F}_q(t)$ at all places of
$\mathbb{F}_q(t)$ in the sense of the classical theory.  The element $\tau \in \mathbb{F}_q^d(\tau)$ satisfies  
 the rule $
\tau \lambda = \lambda^q\tau$ and $\tau^d = t$.  But for all places $p=p(t)
\neq 0, \infty$ we have $v_p(t) = 0$, which implies that the central cyclic
algebra $\mathbb{F}_q^d(\tau)$ is unramified at all places $p \neq 0,
\infty$. On the other hand, computing the invariants at the places $t = 0, \infty$ we obtain for the
invariants
\begin{align*}
&inv_0( \mathbb{F}_q^d(\tau) /  \mathbb{F}_q(t)  = \frac{1}{d}\\
&inv_{\infty}( \mathbb{F}_q^d(\tau) /  \mathbb{F}_q(t)  = -\frac{1}{d}
\end{align*}
We have obtained
\begin{theorem}\label{96}$:$\\
$\mathbb{F}_{q^d}(\tau)$ is up to isomorphism the unique central algebra over 
$\mathbb{F}_q(t)$ with the following properties:\\  $i) \quad
\mathbb{F}_{q^d}(\tau) $ is unramified at all places $p=p(t)
\neq 0, \infty$ of $\mathbb{F}_q(t)$.\\
$ii)$ It has invariants  $\frac{1}{d}$,
$-\frac{1}{d}$ at $t = 0$ resp. $t = \infty$.
\end{theorem}
\begin{proof}
All the properties have been shown, the uniqueness statement is a part of the
theorem of Hasse - Brauer - Noether see 
\cite[chapter XIII,3.Theorem 2 and 6. theorem 4]{we1}
\end{proof}
\begin{lemma} $\mathbb{F}_{q^d}\{\tau\}$ is a maximal $\mathbb{F}_q[t]$-order
 in $\mathbb{F}_{q^d}(\tau)$.
\end{lemma}
\begin{proof} We have to check only the maximality. Assume that $R$ is an order
  with $\mathbb{F}_{q^d}\{\tau\}\subseteq R$. By definition of the concept of order, $R$ is
  finite over $\mathbb{F}_q[t]$, so there exist a nonzero element $f\in
  \mathbb{F}_q[t]$ such that $Rf\subseteq \mathbb{F}_{q^d}\{\tau\}$. As $
  \mathbb{F}_{q^d}\{\tau\}$ is a PID (left and right), there exists an element
  $\mu$ in   $\mathbb{F}_{q^d}\{\tau\}$ with
  $R=\mathbb{F}_{q^d}\{\tau\}\mu$.
  Thus $Rf=\mathbb{F}_{q^d}\{\tau\}\mu f^{-1}$. 
 Now since $R$ is a domain the proof is complete.\end{proof}
Finally we give another description of the skew polynomial ring
  $\mathbb{F}_{q^d}\{\tau\}$.\\ Consider an algebraic closure $\overline{\mathbb{F}}_q$
  of $\mathbb{F}_q$.   $\overline{\mathbb{F}}_q$is a vector space over the field $\mathbb{F}_q$
   and we have the ring of vector space  endomorphisms
  $\text{End}_{\mathbb{F}_q}(\overline{\mathbb{F}}_q)$.
 As mentioned above, there is also the skew polynomial ring
  $\overline{\mathbb{F}}_q\{\tau\}$, satisfying in particular again the
  rule $\tau\lambda = \lambda ^q \tau$ for $\lambda \in \overline{\mathbb{F}}_q$. We
  choose an embedding  \[ \mathbb{F}_{q^d}\{\tau\} \hookrightarrow
  \overline{\mathbb{F}}_q\{\tau\}\]
by choosing a homomorphism $ \mathbb{F}_{q^d} \rightarrow
  \overline{\mathbb{F}}_q$. \\
With any polynomial $f(\tau) = \sum_{i=0}^n a_i\tau^i$, $ f(\tau) \in
  \overline{\mathbb{F}}_q\{\tau\}$ we associate the
  polynomial endomorphism 
\begin{equation}\label{varphi11}
\varphi_f: \overline{\mathbb{F}}_q \longrightarrow
  \overline
{\mathbb{F}}_q,\quad x \mapsto  \sum_{i=0}^n a_i x^{q^i} \end{equation}
\begin{proposition}
$ (i)\qquad \varphi : \bar{\mathbb{F}}_q \{\tau\}\longrightarrow\quad 
  \text{End}_{\mathbb{F}_q}(\overline{\mathbb{F}}_q)$
 \[ f(\tau) = \sum_{i=0}^n
  a_i\tau^i \mapsto \varphi_f\]
 as above, is an injective homomorphism.\\
$ii) \qquad \ker{\varphi_f} = \{ x\in \bar{\mathbb{F}}_q :
  \varphi_f(x) =  \sum_{i=0}^n a_i x^{q^i} = 0 \}$\\ is an $\mathbb{F}_q$-vector
 subspace of $\bar{\mathbb{F}_q }$. One has: \[
 \dim_{\mathbb{F}_q}{\ker{(\varphi_f)}}\leq n \]  
$iii) \qquad \varphi_f$ is injective, iff $f$ is purely inseparable as a
 polynomial, that is, $f(\tau)$ is of the form $ c\tau^n$ with $c \neq 0$.\\
$iv)$ If $f(\tau)\neq 0$, $ \varphi_f$ is surjective.
 \end{proposition}
\begin{proof}
$i)\quad ii)$ and $iii)$ are evident, $iv)$ follows immediately, because the
polynomial equations \[\sum_{i=0}^n a_i x^{q^i} = c \] have solutions for all
$c \in \bar{\mathbb{F}}_q $, iff $ f(\tau) \neq 0$. 
\end{proof}
\begin{proposition}\label{500shah}
 Given a finite dimensional $\mathbb{F}_q$-vector space $V \subset
 \bar{\mathbb{F}}_q $, there is a polynomial $ f(\tau) \in
 \bar{\mathbb{F}}_q\{\tau\}$, unique up  to a scalar from
 ${\bar{\mathbb{F}}_q}^\times$, such that \\
$i) \quad f(\tau)$ is not divisible by $\tau$ $($left or right would be
 equivalent for this$)$.\\
$ii) \quad \ker{(\varphi_f)} = V \subset \bar{\mathbb{F}}_q$,\\
$iii) \quad f(\tau)$ moreover can be chosen to be in $\mathbb{F}_{q^d}\{\tau\}$
 iff $V \subset\bar{\mathbb{F}}_q $ satisfies $\varphi_{\tau^d}(V) = V$, that
 is, the map \[\bar{\mathbb{F}}_q \longrightarrow\quad \bar{\mathbb{F}}_q\] \[
 x \mapsto x^{q^d}  = \varphi_{\tau}(x) \] maps $V$ bijectively onto itself.
\end{proposition}
\begin{proof}
Given  $V \subset\overline{\mathbb{F}}_q $  a finite dimensional vector space
over $\mathbb{F}_q$, we define the polynomial \[ p(x): = \prod_{v\in V}(x - v)
\] By induction with respect to the dimension $\dim_{\mathbb{F}_q}(V)$, it is easy to
show that $p(x)$ is a $ \mathbb{F}_q$-linear $($in particular additive$)$
polynomial function, which therefore is of the form $p(x) = \varphi_f(x)$ for
a skew polynomial  $ f(\tau) \in \bar{\mathbb{F}}_q\{\tau\}$.\\
If there would be another such element  $ g(\tau) \in
 \overline{\mathbb{F}}_q\{\tau\}$ satisfying $\ker{\varphi_g} = V$, then the
 additive polynomial functions $\varphi_f,\quad \varphi_g:\overline{\mathbb{F}}_q
 \longrightarrow \bar{\mathbb{F}}_q$ would have the same zeros in
 $\overline{\mathbb{F}}_q$, furthermore they have the same degree. Therefore there
 exists a $c \in {\overline{\mathbb{F}}_q}^\times$ satisfying $c \varphi_f =
 \varphi_g$.\\
$iii)$ If $V \subset \overline{\mathbb{F}}_q $ satisfies $\varphi_{\tau^d}(V) = V$,
that is, the map $\overline{\mathbb{F}}_q\longrightarrow \overline{\mathbb{F}}_q$, $x
\mapsto x^{q^d} = \varphi_\tau(x)$  maps $V$ bijectively onto itself, then
$\varphi_f$, as constructed above, satisfies 
\begin{align*}
\varphi_{\tau^d}\varphi_f(x) & = {(\prod_{v\in V}(x - v))}^{q^d} \\
& = \prod_{v\in V}(x^{q^d} - v^{q^d})\\
& = \prod_{v\in V}(x^{q^d} - v) \quad (\text{as} \qquad \tau^d(V) = V)\\
& = \varphi_f(\varphi_{\tau^d}(x))
\end{align*}  
 for all $x \in \bar{\mathbb{F}}_q$. This implies immediately  
\begin{equation*}
 \tau^df(\tau) = f(\tau)\tau^d
\end{equation*}
which implies, that $f(\tau)$ is of the form
\begin{equation*}
 f(\tau) = \sum_{i=0}^na_i \tau^i \qquad \text{where} \qquad a_i \in \mathbb{F}_{q^d}.
\end{equation*}
\end{proof} 
\begin{proposition}\label{1010}
If $f(\tau),\quad g(\tau) \in \overline{\mathbb{F}}_q\{\tau\}$ satisfy 
$\ker{(\varphi_f)} = \ker{(\varphi_g)}$ and if   $f(\tau)$ is not divisible by
$\tau$, then there is $ m \in \mathbb{N}$, such that 
 \begin{equation*}
 g(\tau) = \tau^m f(\tau)
\end{equation*}
\end{proposition}
\begin{proof}
Given $g(\tau)$, one can find $m$ maximal, such that $g(\tau) = \tau^{m}
\widetilde{f}(\tau)$ but then $\varphi_{\widetilde{f}}$ is a separable polynomial function
satisfying $\ker{(\varphi_{\widetilde{f}})} = \ker{(\varphi_g)}$ and therefore
also $\ker{(\varphi_{\widetilde{f}})} = \ker{(\varphi_f)}$, furthermore
$\varphi_{\widetilde{f}}$, $\varphi_f$ are separable polynomial functions,
which we can assume to have highest coefficient $1$.  Therefore we
obtain $\widetilde{f} = f$ and therefore also  
\begin{equation*}
 g(\tau) = \tau^m f(\tau)
\end{equation*}
\end{proof}
\begin{proposition}\label{104}
$i)$ Suppose $f_1(\tau),f(\tau) \in \overline{\mathbb{F}}_q\{\tau\}$, are
separable polynomials. Then there is a separable polynomial $f_2(\tau)$
satisfying 
\begin{equation*}
f_2(\tau)f_1(\tau) = f(\tau)\quad\text{iff} \quad \ker{(\varphi_{f_1})}
\subseteq \ker{(\varphi_{f})}\qquad\text{as $\mathbb{F}_q$-subspaces.} 
\end{equation*}
$ii)$ If $f_1(\tau), f(\tau) \in \mathbb{F}_{q^d}\{\tau\}$, are
separable polynomials, then there is \[f_2(\tau)\in \mathbb{F}_{q^d}\{\tau\}\]
satisfying 
\begin{equation*}
f_2(\tau)f_1(\tau) = f(\tau)
\end{equation*}
\end{proposition}
\begin{proof}
 $(\Longrightarrow )$  This direction is trivial. \\
Conversely, suppose that, $\ker{(\varphi_{f_1})} \subseteq \ker{(\varphi_f)}$.
Consider the separable polynomial \[ \varphi_{f_1} :  \overline{\mathbb{F}}_q
\longrightarrow  \overline{\mathbb{F}}_q\] 
and denote for
 \begin{equation*}
 V:= \ker{(\varphi_f)},\qquad \varphi_{f_1}(V) = : \overline{V}\subset
 \overline{\mathbb{F}}_q 
\end{equation*}
Obviously $\bar{V}$ is given as $\mathbb{F}_q$-subspace in $\overline{
\mathbb{F}}_q$ and we can find accordingly  $f_2(\tau) \in
\overline{\mathbb{F}}_q\{\tau\}$, such that $ \varphi_{f_2}$ is separable and
$\ker{\varphi_{f_2}} =  \bar{V}$. Then
 \[ \varphi_{f_2}\varphi_{f_1} : \overline{
  \mathbb{F}}_q \longrightarrow  \overline{\mathbb{F}}_q\] 
has kernel $V$ and is again separable with highest coefficient $1$, Therefore
we obtain $ \varphi_f = \varphi_{f_2}\varphi_{f_1}$, which implies immediately
$f=f_2 f_1$.  This shows $i)$. \\
$ii)$ is an immediate consequence of $i)$.
\end{proof}
\begin{cor}\label{67}
Suppose, $ f(\tau) \in \mathbb{F}_{q^d}\{\tau\}$ is not divisible by $\tau$, that
is, the corresponding polynomial function $ \varphi_f$ is separable.\\
Decompositions of the form 
\begin{equation*}
f(\tau) = f_1(\tau)\ldots f_r(\tau)
\end{equation*}
in $ \mathbb{F}_{q^d}\{\tau\}$ are in bijective correspondence  with $t =
\tau^d$-invariant flags of $\mathbb{F}_q$-linear subspaces,
\begin{equation*}
 0 \subset W_1 \subset \ldots \subset W_r = V 
\end{equation*}
where $\qquad W_{r-j} = \ker{(\varphi_{f_j} \ldots \varphi_{f_r})}\qquad$ for $j= 1,\ldots, r$.
\end{cor}
\begin{proof}
We can assume in the corollary, that $f(\tau)$ and all of the $f_j$ for $j =
1,\ldots, r$ have  highest coefficient $1$.  We then have the map above
associating with a decomposition 
\begin{equation*}
f(\tau) = f_1(\tau)\ldots f_r(\tau)
\end{equation*}
the corresponding flag of subspaces 
\begin{equation*}
 0 \subset W_1 \subset \ldots  \subset W_r = V 
\end{equation*}
where $ V = \ker{\varphi_f}$ and $W_j = \ker{(f_{r -j +1}(\tau))\ldots f_r(\tau))}$.\\
Conversely, suppose that, the flag of $t$-invariant subspaces 
\begin{equation*}
 0 \subset W_1 \subset \ldots \subset W_r = V 
\end{equation*}
is given. We find $f(\tau) \in \overline{\mathbb{F}}_q\{\tau\}$, separable, with highest
  coefficient $1$, such that $\ker{\varphi_f} = V$ holds. 
  $f(\tau)\in \overline{\mathbb{F}}_q\{\tau\}$ is an element in
  $\mathbb{F}_{q^d}\{\tau\}$, because   $V = \ker{\varphi_f}$ is invariant under $t =
  \tau^d$ $($not elementwise however$)$.\\
Similary we find $\widetilde{f}_j(\tau) \in\mathbb{F}_{q^d}\{\tau\}$, such that
  $\ker{(\varphi_{\widetilde{f}_j})} = W_j$, By repeated application of
  Proposition.\ref{104}, we can conclude:
\begin{align*}
 \widetilde{f}_r(\tau) &= f_1(\tau)\ldots f_r(\tau)\\
 \widetilde{f}_{r-1}(\tau) &= f_1(\tau)\ldots f_{r-1}(\tau)\\
&\qquad\vdots\\
 \widetilde{f}_1(\tau) &= f_r(\tau)
\end{align*}
where the $f_j(\tau) \in\mathbb{F}_{q^d}\{\tau\}$.  This shows the Corollary.  
\end{proof}

\section{Arithmetic groups associated to the  division algebra  $\mathbb{F}_{q^d}(\tau)$}
 We consider again the division algebra of skew polynomials $D =
 \mathbb{F}_{q^d}(\tau)$ with center $\mathbb{F}_q(t)$ as in the last
 section. \\ Associated with this algebra are the algebraic groups $D^\times$,
 $D^{(1)}$ and $D^\times/Z$ given  as group functors on the category of
 $\mathbb{F}_q(t)$-algebras $R$ $($that is, there is a homomorphism of
 commutative algebras with unit elements \\
 $\mathbb{F}_q(t)\longrightarrow R$$)$  given by :
\begin{equation}\label{Dtimes}
D^\times(R):= (\mathbb{F}_{q^d}(\tau)\otimes_{\mathbb{F}_q(t)}R)^\times
\end{equation}
the group of units, and similarly
\begin{equation}\label{D1}
D^{(1)}(R):= (\mathbb{F}_{q^d}(\tau)\otimes_{\mathbb{F}_q(t)}R)^{(1)} :=
\{x\in (\mathbb{F}_{q^d}(\tau)\otimes_{\mathbb{F}_q(t)}R)^{(1)} | nr(x)=1\} 
\end{equation}
where $ nr : D^\times(R) \longrightarrow R^\times$ is the reduced norm of the
central simple algebra $\mathbb{F}_{q^d}(\tau)$ over $\mathbb{F}_q(t)$, seen as
a polynomial map and extended by $R$. \\ Finally we have the group functor:
\begin{equation*}
R \mapsto (D^\times/Z)(R)= D^\times(R)/Z(R^\times)
\end{equation*}
\begin{remark}
Other notation for these groups are $\text{\rm{GL}}(1,D)$, $SL(1,D)$ and $PGL(1,D)$. For a
proof, that these group functors are in fact representable by algebraic groups
see \cite{Borel}.
\end{remark}
Besides these algebraic groups, we have the corresponding group schemes
 ${\underline{D}}^\times$, ${\underline{D}}^{(1)}$ and ${\underline{D}}^{\times}/Z$ 
over the ring $\mathbb{F}_q[t]$ respectively. Over its spectrum space
 $\mathbb{F}_q[t]$, given similarly by the group valued functors   :
\begin{equation*}
\mathbb{F}_q[t]-\text{Alg} \qquad \longrightarrow \qquad\text{groups}
\end{equation*}
From the category of commutative $\mathbb{F}_q[t]$-algebras to the category of
groups, given by
\begin{equation*}
R \mapsto D^\times(R)= (D\otimes_{\mathbb{F}_q[t]}R)^\times
\end{equation*}
The arithmetical groups we are considering here can be described as
follows. Suppose $p(t)\in \mathbb{F}_q[t]$ is an irreducible polynomial.
Denote by $\mathcal{O} := \mathbb{F}_q[t][\frac{1}{p(t)}]$ the localization of
the polynomial ring $\mathbb{F}_q[t]$ with respect to the multiplicative
system $\mathcal{S}:= \{p(t)^n | n= 0,\ldots\}$, that is, one considers
the rational functions in $ \mathbb{F}_q(t)$, whose denominators are powers of
$p(t)$. \\ The basic arithmetic groups in our situation are  then
$\underline{D}^\times(\mathcal{O})$,$\underline{D}^{(1)}(\mathcal{O})$,$(\underline{D}^\times/Z)(\mathcal{O})$,
which are given explicitly as 
\begin{align*}
\underline{D}^\times(\mathcal{O}) =&
(\mathbb{F}_q^d\{\tau\}\otimes_{\mathbb{F}_q[t]}\mathcal{O})^\times\\=&(\mathbb{F}_q^d\{\tau\}\otimes_{\mathbb{F}_q[t]}\mathbb{F}_q[t][\frac{1}{p(t)}])^\times
. 
\end{align*}
Furthermore 
\begin{equation}\label{O1}
\underline{D}^{(1)}(\mathcal{O}) = \{ x\in\underline{D}^\times(\mathcal{O})
\quad|\quad nr(x) = 1 \}
\end{equation}
and $\underline{D}^\times/Z)(\mathcal{O}) =
\underline{D}^\times(\mathcal{O})/\mathcal{O}^\times$.\\ 
 We study these groups in the usual way by their operation on the product of the
Bruhat-Tits building of the algebraic group $D^\times (\text{respectively}, D^{(1)},
D^\times/Z)$ at the primes missing, which in this case are  $(p(t))$ and $\infty$ of $
\mathbb{F}_q(t)$.\\
As the division algebra $ \mathbb{F}_{q^d}(\tau)$ is totally ramified
at $\infty$, the corresponding Bruhat-Tits building is just a point. It is
therefore sufficient to consider the   Bruhat-Tits building at the prime
$(p(t))$. This is the building associated to the algebraic group $ D^\times
\otimes_{\mathbb{F}_q(t)} \mathbb{F}_q(t)_{(p(t))}$,  where
$\mathbb{F}_q(t)_{(p(t))}$ is the completion of $\mathbb{F}_q(t)$ at the prime
$p(t): = p$. (and similarly for the groups $D^{(1)}$, $D^\times/Z)$.\\
Because $D$ is unramified at $(p(t))= (p)$, it follows, that, we have an
isomorphism  
\begin{equation*}
D\otimes_{\mathbb{F}_q(t)} \mathbb{F}_q(t)_{p}\quad
\widetilde{\longrightarrow} \quad \mathbb{M}( d ; \mathbb{F}_q(t)_{p}) .
\end{equation*}
Therefore, we have an induced isomorphism 
\begin{equation*}
(D\otimes_{\mathbb{F}_q(t)} \mathbb{F}_q(t)_{p})^\times \quad
\widetilde{\longrightarrow} \quad \text{\rm{GL}}( d ; \mathbb{F}_q(t)_{p}) .
\end{equation*}
We also have an induced embedding 
\begin{align*}
\ \Gamma &=(\mathbb{F}_{q^d}\{\tau\}[\frac{1}{p(t)}])^\times\\
&=
(\mathbb{F}_{q^d}\{\tau\}\otimes_{\mathbb{F}_q[t]}\mathbb{F}_q[t][\frac{1}{p(t)}])^\times\\
&\hookrightarrow \qquad
(\mathbb{F}_{q^d}(\tau)\otimes_{\mathbb{F}_q(t)}\mathbb{F}_q(t)_{p})^\times \quad\cong
  \quad  \text{\rm{GL}}( d ; \mathbb{F}_q(t)_{p})
\end{align*}
and using this, an action of $\Gamma$ and its subgroups on the Bruhat-Tits
building corresponding to $p= p(t)$ respectively also for
 For $\text{\rm{GL}}( d ; \mathbb{F}_q(t)_{p})$ and the related subgroups  $SL( d ;
 \mathbb{F}_q(t)_{p})$ for $D^{(1)}$ and $PGL( d ; \mathbb{F}_q(t)_{p})$ for
 $D^\times/Z$. As described in section \ref{2Baffine}, that is the building $X_.(p) : =
 X_.({\mathbb{F}_q(t)_{p}}^d)$, associated to the vector space
 ${\mathbb{F}_q(t)_{p}}^d$ over the locally compact topological field
 $\mathbb{F}_q(t)_{p}$.\\
The problem we have is to understand the quotient $\Gamma \backslash X_.(p)$ for the
group $\Gamma$ under consideration. To be able to do this we add some further
considerations. \\
First, because $\mathbb{F}_{q^d}\{\tau\}$  is a maximal order over
$\mathbb{F}_q[t]$ in $\mathbb{F}_{q^d}(\tau) = D$, we can choose the isomorphism
  above in such away that it induces an isomorphism of the corresponding
  local orders
\begin{equation*}
\mathbb{F}_{q^d}\{\tau\}\otimes_{\mathbb{F}_q(t)}\mathbb{F}_q[t]_{p}\quad
\widetilde{\longrightarrow}\quad \mathbb{M}( d ; \mathbb{F}_q[t]_{p}),
\end{equation*}
where as before $\mathbb{F}_q[t]_{p}$ is the valuation ring of
$\mathbb{F}_q(t)_{p}$.\\
We denote $L_0 := {\mathbb{F}_q[t]_{p}}^d$, the standard lattice.
For any lattice $ L\subset  {\mathbb{F}_q(t)_{p}}^d$ over
$\mathbb{F}_q[t]_{p}$, we consider  $\text{Hom}_{\mathbb{F}_q[t]_{p}}(L,L_0)$.\\
This is in an obvious way a left module over the ring
\[\text{End}_{\mathbb{F}_q[t]_{p}}(L_0) \cong
\mathbb{F}_{q^d}\{\tau\}\otimes_{\mathbb{F}_q(t)}\mathbb{F}_q[t]_{p} \]
\begin{proposition}
There is a bijective correspondence between lattices \  \\$ L\subset
{\mathbb{F}_q(t)_{p}}^d$ over $\mathbb{F}_q[t]_{p}$ and
$\mathbb{F}_{q^d}\{\tau\}$-left modules $M$ equipped additionally with an
isomorphism \[\phi :
\mathbb{F}_q[t][\frac{1}{p(t)}]\otimes_{\mathbb{F}_q[t]}M\quad
\widetilde{\longrightarrow}\quad \mathbb{F}_{q^d}\{\tau\}[\frac{1}{p(t)}] \]
which are free of rank 1 over $\mathbb{F}_{q^d}\{\tau\}_p = \mathbb{F}_q[t]_{p}\otimes_{\mathbb{F}_q[t]}\mathbb{F}_{q^d}\{\tau\}$
\end{proposition}
\begin{proof}
Given a pair $(M,\phi)$, we have to reconstruct the lattice  $ L\subset
{\mathbb{F}_q(t)_{p}}^d$ over $\mathbb{F}_q[t]_{p}$.  Denoting $\phi_p$ the
obvious extension of $\phi$ to \[ \phi_p :
\mathbb{F}_q(t)_{p}\otimes_{\mathbb{F}_q(t)}M\quad
\widetilde{\longrightarrow}\quad \mathbb{F}_{q^d}(\tau)_p \]
and the $\tilde{\phi_p}$ as the following composition of maps:
\[\mathbb{F}_q(t)_{p}\otimes_{\mathbb{F}_q(t)}M\quad
\widetilde{\longrightarrow}\quad \mathbb{F}_{q^d}(\tau)_p\quad
\widetilde{\longrightarrow}\quad \mathbb{M}(d;\mathbb{F}_q(t)_{p}) \]
where the first map is $\phi_p$, and the second is the isomorphism fixed above. We consider
the restriction map \[ \tilde{\phi_p}_| : M_p =
\mathbb{F}_q[t]_{p}\otimes_{\mathbb{F}_q[t]} M\longrightarrow \mathbb{M}(d;\mathbb{F}_q(t)_p) \]
The image is a free module of rank one over
$\text{End}(\mathbb{F}_q[t]_{p}^d)$. It is immediate to see that there
exists a unique local lattice $ L\subset {\mathbb{F}_q(t)_{p}}^d$, such that
\[ \tilde{\phi_p}_|(M_p) = \text{Hom}_{\mathbb{F}_q[t]_{p}}(L,L_0) \] 
 holds, where $L_0 = {\mathbb{F}_q[t]_{p}}^d$. This is a version of the Morita-
equivalence. Conversely, given the lattice $L$, we obtain an
$\mathbb{F}_{q^d}\{\tau\}$-module $M$ in obvious way from the \[ \mathbb{F}_q[t][\frac{1}{p(t)}]\otimes_{\mathbb{F}_q[t]}M\quad
\widetilde{\longrightarrow}\quad \mathbb{F}_{q^d}\{\tau\}[\frac{1}{p(t)}] \qquad\text{and}\]  
\[
{\mathbb{F}_q[t]}_p\otimes_{\mathbb{F}_q[t]}M\widetilde{\longrightarrow}\quad
\text{Hom}_{\mathbb{F}_q[t]_{p}}(L,L_0) \]
with the canonical identification. It is immediate to see, that these two
constructions are inverse to each other.
\end{proof}
\begin{proposition}\label{bf}
The group $\Gamma = (\mathbb{F}_{q^d}\{\tau\}[\frac{1}{p(t)}])^\times$ acts
transitively on the set of lattices $ L\subset
{\mathbb{F}_q(t)_{p}}^d$ over $\mathbb{F}_q[t]_{p}$.\end{proposition}
\begin{proof} Consider the $\mathbb{F}_{q^d}\{\tau\}$-module $M$ given by the
  pair \[M^{(p)}:= (\mathbb{F}_{q^d}\{\tau\}[\frac{1}{p}],\text{Hom}(L,L_0) )\] 
$($in the sense of the discussion above$)$. Now any such module $($as a left
$\mathbb{F}_{q^d}\{\tau\}$-module $)$ is isomorphic to
$\mathbb{F}_{q^d}\{\tau\}$. Such an isomorphism $\alpha$ induces an isomorphism
$\alpha^{(p)}$ of $\mathbb{F}_{q^d}\{\tau\}_p$-modules 
\begin{equation*} 
M^{(p)}:=
  \mathbb{F}_{q^d}\{\tau\}[\frac{1}{p}] \longrightarrow
  \mathbb{F}_{q^d}\{\tau\}[\frac{1}{p}]
\end{equation*}
Any such  isomorphism is given as right multiplication by a unit \[g \in
{\mathbb{F}_{q^d}\{\tau\}[\frac{1}{p(t)}]}^\times\]
As $\alpha$  induces also an isomorphisms of
${\mathbb{F}_{q^d}\{\tau\}}_p$-modules  
\begin{equation*} 
M_p = \text{Hom}_{{\mathbb{F}_q[t]}_p}(L,L_0) \widetilde{\longrightarrow} \quad
\text{End}_{{\mathbb{F}_q[t]}_p}(L_0),
\end{equation*}
it is immediate to see that this is equivalent to the fact, that \[ g(L) = L_0\]
But this shows the transitivity of the action of  $\Gamma$ on the set of
lattices. 
\end{proof}
\begin{cor}\label{fb}
The group $\underline{D}^\times(\mathcal{O})^\times/Z =
(\mathbb{F}_q\{\tau\}[\frac{1}{p}])^\times/Z$ acts transitively on the set of
vertices $X_0({\mathbb{F}_q(t)_p}^d)$ of the building $X_{\bs{.}}({\mathbb{F}_q(t)_p}^d)$.
\end{cor}
\begin{proof}
We have $Z = (\mathbb{F}_q[t][\frac{1}{p}])^\times$, which acts by scalar
multiplication on the set of lattices. Therefore the group
$(\mathbb{F}_{q^d}\{\tau\}[\frac{1}{p}])^\times/Z$ induces an action on
$X_.({\mathbb{F}_q(t)_p}^d)$, which is transitive on the set of the vertices
$X_0({\mathbb{F}_q(t)_p}^d)$ by proposition \ref{bf}. 
\end{proof}
\begin{cor}\label{finite-simplex} For any subgroup $\Gamma$ of finite index in $\Gamma(1) =
  (\mathbb{F}_{q^d}\{\tau\}[\frac{1}{p}])^\times/Z$ the quotient $\Gamma \backslash
  X_.({\mathbb{F}_q(t)_p}^d)$ is a finite simplicial complex.\end{cor}

\begin{proof} This is an immediate consequence of Corollary \ref{fb}.
\end{proof}
\begin{remark}
In fact this is a very special case of Godement's compactness theorem
mentioned earlier, but in our
situation it can be shown in a direct way as above.\\ We consider in the group
$\Gamma(1) = (\mathbb{F}_{q^d}\{\tau\}[\frac{1}{p}])^\times/Z$  the following
subgroup $\Gamma(\tau)$. We  consider the composition of group homomorphisms 
\begin{equation}\label{3.1.6}
\Gamma(\tau) := \ker{(\Gamma(1) \quad \longrightarrow \quad
  \mathbb{F}_{q^d}\{\{\tau\}\}^\times/Z})  \longrightarrow \quad \mathbb{F}_{q^d}^\times/\mathbb{F}_q^\times
\end{equation} where the first homomorphisms corresponds to the embedding to
  the place $t = 0$ considered easier and the second homomorphism is the
  evaluation homomorphism for $\tau = 0$. $\Gamma(\tau)$ is the kernel of the
  composition of these homomorphisms.
\end{remark}
\begin{proposition}\label{107}
Let $\Gamma(\tau)$ be as \ref{3.1.6}, then\\
$i)$  $\Gamma(\tau)$ is a torsion free group,\\
$ii)$ $\Gamma(\tau)$ acts fixed point free on the simplicial complex
$X_.({\mathbb{F}_q(t)_p}^d)$ and also on its realization
$|X_.({\mathbb{F}_q(t)_p}^d)|$. In particular no simplex is mapped to itself
in a nontrivial way  under the action of $\Gamma(\tau)$.
\end{proposition}
\begin{proof}
$i)$ We consider the embedding  
\[(\mathbb{F}_{q^d}\{\tau\}[\frac{1}{p}])^\times/Z \hookrightarrow  \mathbb{F}_{q^d}\{\{\tau\}\}^\times/Z \]
Any element $ g \in \Gamma(\tau)$ $g \neq 1$ is mapped to a power series in
$\tau$ of the form $(1 + a_i\tau^i + \text{higher terms in}\quad \tau)$, where $a_i
\in \mathbb{F}_{q^d}$, $a_i \neq 0$. It is immediate to see that  $(1 +
a_i\tau^i + \text{higher terms in} \tau)^n \neq 1$ if $(n,Char(\mathbb{F}_q )
= 1$. On the other hand for any $j \geq 1$ 
\begin{align*}
(1 +a_i\tau^i + \text{higher terms in}\quad \tau)^{p^j} &= (1 + a_i\tau^i (1 +
b_1\tau + \ldots))^{p^j}\\& = 1 + (a_i)^{p^j} +\ldots
\end{align*} 
Thus the lowest term in $\tau$ is $(a_i)^{p^j}$, which obviously is not
zero. This shows $i)$.\\
$ii)$ If $ g \in \Gamma(\tau)$ $g \neq 1$, stabilizes a simplex, it would have
a fixed point in the realization $X_.({\mathbb{F}_q(t)_p}^d)$. As the
stabilizer is discrete $($as a subgroup of $\Gamma(\tau))$ and compact $($
being a closed subset in the compact stabilizer of a point of the
building $X_.({\mathbb{F}_q(t)_p}^d))$, it follows, that such a stabilizer
group is finite. This implies immediately that $g$ has a finite order, which is a
contradiction.\end{proof}
\begin{proposition} $\Gamma(\tau)$  is transitive on the set of the vertices of
  the building $X_.({\mathbb{F}_q(t)_p}^d))$.\end{proposition}
\begin{proof} We have seen above, that the group $\Gamma =
  (\mathbb{F}_q\{\tau\}[\frac{1}{p}])^\times$ is transitive on the set of
  vertices $X_0({\mathbb{F}_q(t)_p}^d))$.  Therfore, given a vertex $<L> \in
  X_0({\mathbb{F}_q(t)_p}^d))$,  we find  $ g \in \Gamma(\tau)$, $ gL = L_0$, where $L_0 = {\mathbb{F}_q[t]_p}^d$. If  there is another  $ g' \in
  \Gamma(\tau)$ with $ g'L = L_0$, then we have $g^{-1}g'L_0 = L_0$. It
  follows immediately from Proposition \ref{107} $(i)$ that $g^{-1}g' = 1$ and
  so $ g = g'$. This completes the proof of the corollary. 
\end{proof}
 \begin{cor}\label{e987} $\Gamma(\tau)$ acts simply transitive on the
  building $X_{\bs{.}}({\mathbb{F}_q(t)_p}^d)$.\end{cor}
\begin{proof}
Immediately obtained from the above propositions.\end{proof}
\section{Explicit description of some arithmetic quotients}
In $4.3$ we have introduced the arithmetic group
$\Gamma(\tau)$
\begin{equation*}
\Gamma(\tau) \subset \Gamma(1) =
\biggl(\mathbb{F}_{q^d}\{\tau\}[\frac{1}{p}]\biggr)^\times/Z
\end{equation*}
 which acts simply transitive on the building $X_{\bs{.}}(\mathbb{F}_q(t)_p^d)$ and
 also on its topological realization $|X_{\bs{.}}(\mathbb{F}_q(t)_p^d)|$.\\
We consider now arbitrary linear polynomials $p(t) = (t-\lambda)$ for $\lambda
 \in \mathbb{F}_q$, $\lambda \neq 0$. Associated to the polynomial $p(t) = t -
 \lambda$ is the $\mathbb{F}_q$-linear homomorphism \begin{equation*}
\varphi_p : \overline{\mathbb{F}}_q \longrightarrow \overline{\mathbb{F}}_q,
 \quad x \mapsto (x^{q^d} - \lambda x)
\end{equation*} Denote  $V : = \ker{(\varphi_p)} = \{x
 \in\overline{\mathbb{F}}_q \quad |  \quad x^{q^d} - \lambda x = 0\}$\\
Of course, $V = V(\lambda)$  is explicitly given as $\biggl(\mathbb{F}_{q^d}
\bs{.} \sqrt[q^d]{\lambda}\biggr)$, a one-dimensional
$\mathbb{F}_{q^d}$-sub-vector space in $\overline{\mathbb{F}}_q$.\\
In particular for the case $\lambda = 1$ we have $V(\lambda = 1)
=\mathbb{F}_{q^d} \subset  \overline{\mathbb{F}}_q$. We apply Corollary
\ref{67} from $4.2.$ to this situation :\begin{proposition}
 There is a bijective correspondence between decompositions 
\begin{equation*}
p(t) = f_1(\tau) \ldots f_d(\tau)
\end{equation*}
into $\tau$-linear factors $f_j(\tau) \in \mathbb{F}_{q^d}\{\tau\}$, $j = 1,
\ldots , d$, and arbitrary full flags
\begin{equation*}
 0 \subset W_1 \subset \ldots \subset W_d = V \quad \end{equation*} 
of $\mathbb{F}_q$-linear subspaces $W_j \subset V$ $($such that $\dim{(W_j)} =
j$. This correspondence  is given as in corollary \ref{67}, section $4.3.$
\end{proposition}
\begin{proof} We only have to check that an arbitrary full flag
 \begin{equation*}
 0 \subset W_1 \subset \ldots \subset W_d = V 
\end{equation*} 
is $t = \tau^d$-invariant. But on $V$ the relation $\tau^d = \lambda$
resp. $x^{q^d} = \lambda x$ holds. Therefore the action $x \mapsto x^{q^d}$ on
$V$ is exactly the homothety $V \longrightarrow V,\quad x \mapsto \lambda x$
with $\lambda \in \mathbb{F}_q$. Because $ \lambda \bs{.} W_j = W_j$, as the
$W_j$ are $ \mathbb{F}_q$-vector-spaces in $V$, the result follows.
\end{proof}
We fix now again the standard lattice $L_0 : = (\mathbb{F}_{q^d}\{\tau\})^d
\subset \mathbb{F}_{q^d}(\tau)_p^d$ as well as the associated vertex
$<L_0>$. By Proposition \ref{link0} from Chapter $2$ we have the isomorphism of simplicial complexes
\begin{equation*} 
lk\biggl(<L_0>; X_{\bs{.}}(\mathbb{F}_q(t)_p^d)\biggr)
\end{equation*}
and the Tits building associated to the $\mathbb{F}_q$-vector-space $(L_0/ \pi
L_0)$.  We remind the reader again of the following situation. We have fixed an
isomorphism   
\begin{equation*}
\mathbb{F}_{q^d}\{\tau\}_{(p)} \widetilde{\longrightarrow} \mathbb{M}(d; \mathbb{F}_q[t]_p)
\end{equation*}
of the corresponding completions at $p = p(t)$. Upon doing this we get a
corresponding standard lattice $\mathbb{F}_q[t]_p^d = L_0$, which is acted
upon by $\mathbb{F}_{q^d}\{\tau\}_{(p)}$ resp. $\mathbb{M}(d;
\mathbb{F}_q[t]_p)$ by the standard action. We have the quotient  $(L_0/ \pi
L_0)$, where $\pi = p(t)$, isomorphic to ${(\mathbb{F}_q[t]/(p(t))}^d$ as
$\mathbb{F}_q[t]/(p(t) \cong \mathbb{F}_q$-vector space.\\
Furthermore we have $V = V(p(t)) = \ker{(\varphi_p)} \subset
\overline{\mathbb{F}}_q$, which is also an $\mathbb{F}_q$-vector space of
dimension $d$. Both  $(L_0/ \pi L_0)$ and $V$ are naturally simple modules over
\begin{equation*}
\mathbb{F}_{q^d}\{\tau\}/p(t)\mathbb{F}_{q^d}\{\tau\} \cong \mathbb{M}(d; \mathbb{F}_q[t]/(p(t)).
\end{equation*}  
An isomorphism between these two modules, so compatible with the action
$\mathbb{F}_{q^d}\{\tau\}/(p(t))$ on both sides, by Schur's lemma will be
unique up to a central nonzero element of $\mathbb{F}_q[t]/(p(t) \cong
\mathbb{F}_q$, that is, up to a scalar, $\neq 0$. This implies in particular a
unique isomorphism between the simplicial complexes  \begin{equation*} 
lk\biggl(<L_0>; X_{\bs{.}}(\mathbb{F}_q(t)_p^d)\biggr)
\end{equation*} and the Tits building of the $\mathbb{F}_q$-vector space $V = V(p(t)) = \ker{(\varphi_p)} \subset
\overline{\mathbb{F}}_q$, compatible with the action of
$\biggl(\mathbb{F}_{q^d}\{\tau\}/(p(t))\biggr)^\times$.\\ In particular, to
any neighboring vertex $<L> \in X_0(\mathbb{F}_q[t]_p^d)$, such that
$<L_0,L_1> \in X_1(\mathbb{F}_q[t]_p^d)$ $(1$-simplices$)$, there is, unique
element $\gamma_L \in \Gamma(\tau)$, because $\Gamma(\tau)$ acts simply
transitive on the set of vertices.\\ We can also obtain explicit
descriptions of the elements $\gamma_L$ using the canonical isomorphisms
above.\\ So, suppose we have a   neighboring vertex $<L>$ of  $<L_0>$, so we
assume the situation  \begin{equation*} 
L_0 \supsetneqq L \supsetneqq p(t)L_0.
\end{equation*}
We have $ \gamma_L <L_0> = <L>$ upon multiplying by an appropriate power of
$p(t)$, we can even assume, that $ \gamma_L(L_0) = L$.\\ In particular, we
obtain then  $ \gamma_L(L_0) \subset L_0$ and therefore, because $\gamma_L$ is
free of all other primes, we obtain $\gamma_L \in \mathbb{F}_{q^d}\{\tau\}$. But
then we obtain the induced equality \begin{equation*} 
\gamma_L \bs{.}\biggl(L_0/p(t)L_0\biggr) = \biggl(L/p(t)L_0\biggr),
\end{equation*}
That is, in term of the vector space $V = V(p(t)) \subset
\overline{\mathbb{F}}_q$ : $\gamma_L \bs{.}V = W(\overline{L}) \subset V$.\\
The divisor of $p(t)$ corresponding to $\gamma_L$ will be therefore given by
the equality above. We have proved 
\begin{theorem}
There is canonical bijection between the sets 
\begin{equation*} 
\biggl\{ \gamma_{L} \quad | \quad  \gamma_{L} \in \mathbb{F}_{q^d}\{\tau\},
  \gamma_{L}(L_0) = L, \quad \text{such that} <L_0, L> \in
  X_1(\mathbb{F}_q[t]_p^d)\biggr\}
\end{equation*} and  
\begin{equation*}
\biggl\{ f(\tau) \in \mathbb{F}_{q^d}\{\tau\} \quad | \quad f(\tau) \quad \text{is a
nontrivial divisor of}\quad p(t)\biggr\},
\end{equation*}
\[ \gamma_L \mapsto f(\tau), \]
given by the equation \[ \gamma_L(L_0) / p(t) L_0 \cong f(\tau)\bs{.}V \]
with respect to the canonical identification \[ L_0 / p(t) L_0\quad 
\widetilde{\longrightarrow} \quad V = \ker{(\varphi_p)} . \]
\end{theorem}
 We consider now the Cayley graph of $\Gamma(\tau)$ with respect to the set of 
 generators  
\begin{equation*} 
\biggl\{ \gamma_{L} \quad | \quad  \gamma_{L} \in \mathbb{F}_{q^d}\{\tau\},
  \gamma_{L}(L_0) = L, \quad \text{such that} <L_0, L> \in
  X_1(\mathbb{F}_q[t]_p^d)\biggr\},
\end{equation*}
We denote $\text{Cayley}(\Gamma(\tau); \{\gamma_L\})$ the associated graph.
\begin{theorem}\label{oben0}
There is a canonical isomorphism of graphs 
\begin{equation*}
 \text{Cayley}(\Gamma(\tau); \{\gamma_L\}) \quad \widetilde{\longrightarrow} 
\quad \tau_{\leqslant 1}\biggl(X_{\bs{.}}(\mathbb{F}_q(t)_p^d)\biggr),
\end{equation*}
where $\tau_{\leqslant 1}\biggl(X_{\bs{.}}(\mathbb{F}_q(t)_p^d)\biggr)$ is the
graph underlying the simplicial complex
$X_{\bs{.}}(\mathbb{F}_q(t)_p^d)$.\end{theorem}
This isomorphism is given by the map 
\begin{equation*}
\text{Cayley}(\Gamma(\tau); \{\gamma_L\})_{(0)} = \Gamma(\tau) \quad
\longrightarrow \quad X_0(\mathbb{F}_q(t)_p^d),
\end{equation*}
\begin{equation*}
\gamma \mapsto \gamma <L_0>.
\end{equation*}
\begin{proof} The map above on the $0$-level is a bijection, because
  $\Gamma(\tau)$ is simply transitive on $X_0(\mathbb{F}_q(t)_p^d)$. Two
  vertices $<\gamma>$ and  $<\gamma^\prime>$, where $\gamma, \gamma^\prime \in
  \Gamma(\tau)$ define a $1$-simplex $<\gamma,\gamma^\prime>$ in the Cayley
  graph $\text{Cayley}(\Gamma(\tau); \{\gamma_L\})$, iff there is $\gamma_L
  \in \Gamma(\tau)$, satisfying $\gamma \gamma_L =  \gamma^\prime$. But then 
 \begin{equation*}
<\gamma \gamma_L(L_0), \gamma L_0> = \gamma <\gamma_L(L_0),  L_0>,
\end{equation*}
 which obviously is a one-simplex in $X_{\bs{.}}(\mathbb{F}_q(t)_p^d)$, because $
 <\gamma_L(L_0),  L_0>$ is a one-simplex in
 $X_{\bs{.}}(\mathbb{F}_q(t)_p^d)$ by definition.Furthermore, this isomorphism
 is label-preserving and induces an isomorphism of the corresponding
 hypergraphs resp. simplicial complexes. The converse is equally clear,
 furthermore it is clear, that our isomorphism will  preserve labels, which
 can be introduced via the set of generating elements $\{\gamma_L\}$. By
 proposition \ref{24} from  Chapter $2$ the graph structure in our situation induces in a unique way a
 simplicial structure resp. the structure of a hypergraph. Obviously, our
 isomorphism induces an isomorphism of this structure.\end{proof}
  \begin{theorem}
Suppose, $\Gamma \subset \Gamma(\tau)$ is a normal subgroup of finite index,
$\overline{\gamma_L} := proj(\gamma_L)$, where \[proj : \Gamma(\tau) \quad
\longrightarrow \quad \Gamma(\tau) / \Gamma \] is the canonical
projection. Then the isomorphism of Theorem \ref{oben0} induces an isomorphism of the graph
\begin{equation*}
\text{Cayley}(\Gamma(\tau)/\Gamma; \{\overline{\gamma_L}\}) =  \quad
\longrightarrow \quad \tau_{\leqslant 1}\biggl(\Gamma  \backslash X_{\bs{.}}(\mathbb{F}_q(t)_p^d)\biggr)
\end{equation*}
 This isomorphism of graphs is again label-preserving and induces an
 isomorphism of the corresponding hypergraphs resp. simplicial complexes.
\end{theorem}
\begin{proof}
This follows immediately from Theorem  \ref{oben0}.
\end{proof}
\section{Explicit construction}
In this section we fix $p(t) = 1 - t$  where as before $t=\tau^d$. Our goal is
 a suitable linear factorization of $p(t)$ in
$\mathbb{F}_{q^d}\{\tau\}$. More precisely we hvae: 
\begin{theorem}\label{above108}
 There are $x_1,\ldots, x_{d-1},x_d \in \mathbb{F}_{q^d}$ $($provided $x_d = 1
 )$ such that 
$(a)$ For any $i \in \{0,\ldots,d-1\}$ there exist the linearly indepndent set
 $\{y_{d-i}^{(j)}\}$  $($provided $y_d^{(0)} = 1)$
with \[ \ker{\varphi_{\prod_{j=0}^i(1-x_{d-j}^{1-q}\tau)}}
=\bigoplus_{j=0}^i \mathbb{F}_q \cdot y_{d-i}^{(j)} .\]
$(b)$ The following linear factorization of $p(t)$ in
$\mathbb{F}_{q^d}\{\tau\}$ holds:
\begin{equation} \label{f00}
p(t) = 1 - t = (1 - \tau)(1 - x_{d-1}^{1-q}\tau) \ldots (1- x_1^{1-q} \tau).
\end{equation}
Recall that the $\mathbb{F}_q$-linear map $\varphi$ is defined by \ref{varphi11}.
\end{theorem}
 \begin{proof} Let $i =0$  then $\ker{\varphi_{(1 - \tau)}} = \mathbb{F}_q
 \cdot 1$ so  $y_d^{(0)} = 1$. Assume $($Inductive assumption$)$ now,
 $\{y_{d-i}^{(j)}\}$ is defined for $1 \le j \le i < k$ such that $(1)$
 holds.  We must determine $x_{d-k}$ and
 $\ker{\varphi_{\prod_{i=0}^k(1-x_{d-j}^{1-q}\tau)}}$.
\[ z \in \ker{\varphi_{\prod_{i=0}^k(1-x_{d-j}^{1-q}\tau)}} \Longleftrightarrow
 \varphi_{(1 - x_{d-k}^{1-q}\tau)} \in \ker{\varphi_{\prod_{i=0}^{k - 1}(1-x_{d-j}^{1-q}\tau)}}\] by inductive assumption
\[  \Longleftrightarrow z - x_{d - k}^{1-q} z^q \in \bigoplus_{j=0}^i
\mathbb{F}_q \cdot y_{d-i}^{(j)} \ . \]
Let \[ V_k^\bot := \{ u \in \mathbb{F}_q \ | \
\text{Tr}_{\mathbb{F}_{q^d}/\mathbb{F}_q}(\frac{y_{d-i}^{(j)}}{u}) = 0, \ 0
\le j \le k-1 \} \]
and choose a nonzero element $x_{d-k} \ \in  \ V_k^\bot$. We have for any $ i \le
k$
\[  \varphi_{1-x_{d - k}^{1-q}\tau}(z) = z - x_{d - k}^{1-q} z^q =  y_{d-i}^{(j)} \Longleftrightarrow
(\frac{z}{x_{d-k}}) - (\frac{z}{x_{d-k}})^q = (\frac{y_{d-i}^{(j)}}{x_{d-k}})^q \
      . \]  By Theorem Hilbert $90$, $($see \cite[ page 215 ]{Lang}$)$, we can
      find for $0 \le j \le k$ $\theta_{d-i}^{(j)}$, $($provided
      $\theta_{d-k}^{(0)} = 1)$
   such that
 \[ \theta_{d-i}^{(j)} - (\theta_{d-i}^{(j)})^q =
      (\frac{y_{d-i}^{(j)}}{x_{d-k}})^q , \]
let $y_{d-k}^{(j)} := x_{d-k} \theta_{d-k}^{(j)}$. Then $y_{d-k}^{(j)} \in
\ker{\varphi_{\prod_{i=0}^{k - 1}(1-x_{d-j}^{1-q}\tau)}}$. \\
Cliam : The set $\{ y_{d-k}^{(0)}, y_{d-k}^{(1)}, \ldots, y_{d-k}^{(k)}\}$ is
linearly independent.
\begin{proof}Of cliam: 
\begin{align*}
&\text{Let} \quad  \sum_{j=0}^k c_j y_{d-k}^{(j)} = 0\\
&\Longrightarrow \quad  \sum_{j=0}^k c_j \frac{y_{d-k}^{(j)}}{x_{d-k}} =
    0\\
&\Longrightarrow \quad  \sum_{j=0}^k c_j \theta_{d-k}^{(j)} = 0\ \Longrightarrow \   (\sum_{j=0}^k c_j
\theta_{d-k}^{(j)})^q = 0 \Longrightarrow \  \sum_{j=0}^k c_j (\theta_{d-k}^{(j)})^q = 0\\
&\Longrightarrow \quad  \sum_{j=0}^k c_j (\theta_{d-k}^{(j)}-
(\theta_{d-k}^{(j)})^q) = 0\\
&\Longrightarrow \quad \sum_{j=1}^k c_j \frac{y_{d-k+1}^{(j-1)}}{x_{d-k}} = 0\ \Longrightarrow \quad \sum_{j=1}^k c_j y_{d-k+1}^{(j-1)} = 0
\end{align*}
Now by inductive assumption we know that the set $\{ y_{d-k+1}^{(0)},
y_{d-k+1}^{(1)}, \ldots, y_{d-k+1}^{(k-1)}\}$ is
a linearly independent
 set over $\mathbb{F}_q$. So \[ c_1 = \ldots c_k = 0
\Longrightarrow c_0 = 0.\] 
\end{proof}
We have by definition  $ \dim{V_j^\bot} = d - j$. So we find $x_d = 1, x_{d-1},
\ldots , x_1$ and a linearly indepndence set $\{y_1^{(0)} = x_1, y_1^{(1)},
\ldots ,y_1^{(d-1)}\}$ such that :
\[ \ker{\varphi_{\prod_{j=0}^{d-1}(1-x_{d-j}^{1-q}\tau)}}
=\bigoplus_{j=0}^{d-1} \mathbb{F}_q \cdot y_1{(j)} .\] This proves $(1)$.\\
 In order to prove $(2)$,
 we see that the map $\varphi$   on the both side of \ref{f00}
has the same kernel $($in  this case, the kernel is 
$\mathbb{F}_{q^d} )$. \\Clearly $\ker{\varphi_{1-t}} = \mathbb{F}_{q^d}$,
  and also by ${(1)}$ we have :
\[
\ker{\varphi_{(1 - \tau)(1-x_{d-1}^{1-q}\tau)\ldots(1-x_1^{1-q}\tau)}}
      = \mathbb{F}_{q^d}.\]
 
Thus by Proposition \ref{1010}
there is a non-negative integer $m$ such that
\begin{equation}\label{400t}
 1-t = \tau^m(1 + \sum^{d-1}_{i=1}x_i^{1-q}\tau)(1-x_{d-1}^{1-q}\tau)\ldots(1-x_1^{1-q}\tau)\cdot
\end{equation}
 From the following lemma,  all factors in the left side
of $(\ref{f00})$ have reduced norm equal to $1-t$.  Take the reduced norm from
 both sides of $(\ref{400t})$. Then we must have  $t^m = 1 $, so  $m = 0 $. Thus the proof of
$(2)$ is complete.
\end{proof}
\begin{lemma}\label{4800} For any $x, \in  \mathbb{F}_{q^d}^\times$ we have:
  \[ \quad \text{rn}(1-x^{1-q}\tau) = 1 - t. \]
where rn is the reduced norm defined by \ref{reduced1.1} in Chapter $1$, Section $4$. \end{lemma}
\begin{proof}
We have 
\begin{align*}\label{4900}
\text{rn}(1-x^{1-q}\tau) &=\det{\begin{pmatrix}  
1 &-x^{1-q}  &  & \ldots & 0\\
                 0 & 1 & -x^{q-q^2} & \ldots &
                \\\vdots & \vdots &
                 \vdots & \ddots & \vdots \\ 0 & 0 & 0 & \ldots &-x^{q^{d-2}-q^{d-1}} \\-tx^{q^{d-1}-q^d} &
                 0 & 0 & \ldots & 1
\end{pmatrix}}\\
&= 1 + {(-1)}^{d+1} {(-1)}^d x^{\sum_{\scriptstyle i=0}^{d-1}(q^i - q^{i+1})}t\\
&= 1 - t.
\end{align*}   
\end{proof}
\newpage
\begin{definition}
For any $i = 0 ,\ldots, d-1$ let $F_i$ be the set of all 
\[ \prod_{j=0}^i(1-x_{d-j}^{1-q}\tau) ,\]
where $x_d = 1$ and $x_j \in \mathbb{F}_{q^d}$ $j = d-i,\ldots,d-1$
such that there exists $(x_{d-i+1},\ldots,x_1)$ with
\[  (x_1,\ldots, x_i,x_{i+1},\ldots,x_{d-1}) \in \mathcal{B}^{1 - t}_{q,d}\cdot \]
We define the fundamental domain of $p(t) = 1-t$ as the disjoint union of
$F_i$'s and we denote it by $\text{FUND}_{1-t}$. 
\end{definition}
\begin{cor} There is a bijective correspondence between $F_i'$s in the
  above definition and $Gr_i(\mathbb{F}_q)$ $($the Grassmanian of$\mathbb{F}_{q^d}$ over $\mathbb{F}_q$. 
\end{cor}
\begin{proof}
The discussion in the Proof of Theorem \ref{above108} makes this bijective correspondence clear.\end{proof}
 Now we return to $\Gamma(\tau)$ as \ref{3.1.6}. \\
\begin{definition}\label{Gamma}
If in  $(\ref{3.1.6})$  $p(t) = 1-t$, then  we denote 
instead $\Gamma(\tau)$, $\Gamma^{1-t}(\tau)$, that is: 
\begin{equation*}
\Gamma^{1-t}(\tau) := \ker{\biggl(\Gamma^{1-t}(1) \quad \longrightarrow \quad
  \mathbb{F}_q\{\{\tau\}\}^\times/Z}\biggr) 
\end{equation*}
where
 \begin{equation}\label{gamma1}
   \Gamma^{1-t}(1):= (\mathbb{F}_{q^d}\{\tau\}[\frac{1}{1-t}])^\times /(Z =
   \mathbb{F}_{q^d}^\times)
 \end{equation}
\end{definition}
\begin{remark} From now on, we shall be  working with groups modulo their
   centers. However, our calculations are made in the groups themselves (by
   taking 
 arbitrary liftings).\end{remark} 
\begin{proposition}\label{120p}
 Let $\Gamma^{1-t}(\tau)$ as definition $($\ref{Gamma}$)$. Then we have:\\
 $i)$  $\Gamma^{1-t}(\tau)$ is a torsion free group.\\
 $ii)$ $\Gamma^{1-t}(\tau)$ is a finitely generated group which is generated
 by the set$\text{FUND}_{1-t}$. \\
 $iii)$ The Cayley graph of $\Gamma^{1-t}(\tau)$ with respect to the generator
 set $\text{FUND}_{1-t}$ is isomorphic with the vertex set of Building
 $X_.(\mathbb{F}_q(t)_{1-t}^d)$, i.e.$X_0(\mathbb{F}_q(t)_{1-t}^d \cong PGL(d, \mathbb{F}_q(t)_{1-t}/ PGL(d, \mathbb{F}_q[t]_{1-t})$.
\end{proposition}
\begin{proof}
We can see directly that any element of $\Gamma^{1-t}(\tau)$ can be written as
$1 + \tau^s\mu$ in which $\mu$ is not divisible by $\tau $. Thus the group
$\Gamma^{1-t}(\tau)$ has torsion elements iff there are positive integers $N$
and $m$ respectively such that:
\begin{equation}\label{torsion}
 (1 + \tau^s\mu)^N = (1 - t)^m \cdot
\end{equation}
 Thus we have :
\begin{equation*}
1 + N\tau^s\mu + \tau^s\mu\tau^s\mu + \ldots = 1 - mt + m(m-1)/2t^2 - \ldots  \cdot
\end{equation*}
Now  if  $s\neq d$ we obtain a contradiction, since $\mu$ is not divisible by $\tau $ and
neither is the coefficient of $t$ in the right-hand side of above equation.\\
So the only possibility is that $s = d$. That is $(1 + \tau^s\mu)$ must be in
\[ \Gamma^{(1-t)}(t): = \{ \alpha \in \Gamma^{1-t}(\tau) \quad | \quad \alpha
\equiv 1 \mod{t}\}\] which is by \cite[Chapter $5$]{ddms4}, 
 an analytic torsion
free group. This proves $i)$.  
We have $ \mu \in  \Gamma^{1-t}(\tau) \quad \Longleftrightarrow \quad
rn(\mu) = (1-t)^m$ for some positive integer $m$. Suppose that $s$ is the  maximal
power of $(1-t)$ dividing $\mu$; then $ \gamma = \mu/(1-t)^s$ is not divisible
by $1-t$. But $nr(\gamma) = (1-t)^{m - ds}$,  and on the other hand, since
$\mathbb{F}_{q^d}\{\tau\}$ is a $($left$)$ PID, there exists an element $\theta \in
\mathbb{F}_{q^d}\{\tau\}$ which generates the ideal $(\gamma, 1-t)$. So 
$\gamma = \alpha\theta$,  and using associativity $($ multiplication by a unit
on the left $)$  $\theta$ can be chosen uniquely. Also since $\theta$ divides 
$1-t$ , $nr(\theta)$ divides $(1-t)^m$, which shows that $rn(\theta) = (1-t)^i$
for some $ i=1,\ldots, m $, and $\varphi_{\theta}$  must have the same  kernel as some
$\varphi_{\prod_{j=0}^i(1 - x_{d-j}^{1-q}\tau)}$  $($again since $\theta$
divide $1-t$ and using Proposition \ref{500shah} $iii))$. This shows that by  Proposition
\ref{1010} there exists a non-negative integer $n$ such that $ \theta = \tau^n\prod_{j=0}^i(1 - x_{d-j}^{1-q}\tau)$. Take
reduced norm of both sides to see $ n$ must be $0$. Now we can repeat this discussion with
$\gamma/\theta =\alpha$. This completes the proof of  $ii)$.\\
  $iii)$ is obtained immediately from Corollary \ref{e987}.
\end{proof}
\begin{definition} Let $f\in \mathbb{F}_q[t]$ be an irreducible element, different
  from $t$ and $1-t$. Then we define:
 \begin{equation*}
 \Gamma_f^{1-t}(\tau) := \{\mu \in \Gamma^{1-t}(\tau) \quad | \quad \mu \equiv 1 \mod{f}\}
\end{equation*}
\end{definition}
We are now ready to present our main results of this Chapter.
 \begin{theorem} \label{Hyper} 
{$(\text{Main Explicit Theorem})$} \\
 Let $f\in \mathbb{F}_q[t]$ be an irreducible element, different
  from $t$ and $1-t$. Then the Cayley graph of
$\Gamma^{1-t}(\tau)/\Gamma_f^{1-t}(\tau) $ with respect to the
$\overline{\text{FUND}}_{1-t}$   \\ $($the image of
  $\text{FUND}_{1-t}$ in the quotient group
  $\Gamma^{1-t}(\tau)/\Gamma_f^{1-t}(\tau) )$, is a Ramanujan $(n_1,n_2,\ldots,n_{d-1})$-regular hypergraph, 
where $n_i$ is the number of all $i$-dimensional sub-vector spaces of
$\mathbb{F}_{q^d}$  with the bound $(c_1, \ldots, c_{d-1})$ where
  $c_i=\binom{d}{i}q^{(d-1)/2}$ for $i=1,2,\ldots,d-1$. 
 \end{theorem}
\begin{proof}
We can consider the quotient group \[\Gamma^{1-t}(\tau)/\Gamma_f^{1-t}(\tau) \] as
$\Gamma_f^{1-t}(\tau)\backslash \Gamma^{1-t}(\tau)$.  By above Proposition \ref{120p} $iii)$
 $\Gamma^{1-t}(\tau)$  can be identified with $X_0(\mathbb{F}_q(t)_{1-t}^d)$,
so the hyper graph properties come from the building structure on $X_.(\mathbb{F}_q(t)_{1-t}^d)$, and finiteness from 
Corollary \ref{finite-simplex}. The Ramanujan property is an immediate
consequence of Theorem \ref{abstract} if d is a prime number, and Corollary \ref{cor47}  in the  Chapter $3$.
\end{proof}
$Notation:$ \\
We denote the hypergraph in Theorem \ref{Hyper} with
$Hyp_f(1-t)$. Our goal is, to give a very simple form of $Hyp_f(1-t)$.\\ 
In order to the simplify calculations, we assume that $\deg{f}=dn$ for some
positive integer $n$. We define now the isomorphism :
\begin{equation}\label{pepsi}
\psi:\Gamma_{1-t}(\tau)\longrightarrow
PGL(d, \mathbb{F}_{q^{\deg(f)}})
\end{equation} by
  \begin{equation*}
\psi([\sum_{i=0}^{d-1}\alpha_i\tau^i])=\begin{pmatrix}  
                \overline{\alpha_0} &\overline {\alpha_1} & \overline{\alpha_2} & \ldots & \overline{\alpha_{d-1}}\\
                 \overline{t}\overline{\sigma(\alpha_{d-1})} & \overline{\sigma(\alpha_0)} & \overline{\sigma(\alpha_1)} &
                 \ldots & \overline{\sigma(\alpha_{d-2})} \\\overline{t}\overline{\sigma^2(\alpha_{d-1})} &
                \overline{t}\overline{\sigma^2(\alpha_0)} & \overline{\sigma^2(\alpha_1)} & \ldots &
                 \overline{\sigma^2(\alpha_{d-2})}\\\vdots & \vdots &
                 \vdots & \ddots & \vdots \\ \overline{t}\overline{\sigma^{d-1}(\alpha_1)} &
                 \overline{t}\overline{\sigma^{d-1}(\alpha_2)} & \overline{t}\overline{\sigma^{d-1}(\alpha_3)} & \ldots &
                 \overline{\sigma^{d-1}(\alpha_0)}
\end{pmatrix}
\end{equation*}
where $[\sum_{i=0}^{d-1}\alpha_i\tau^{i-1}]$ is the image
$\sum_{i=0}^{d-1}\alpha_i\tau^{i-1}$ in $\Gamma_{1-t}(\tau)$, and
\\$\quad^{\overline{\quad}}:\mathbb{F}_q[t]\longrightarrow
\mathbb{F}_q[t]/(f)\cong \mathbb{F}_{q^{\deg(f)}}$ is the natural map $($which extended simply over constant field $)$.
We have 
\begin{equation*}
 \ker{\psi} = \Gamma_f^{1-t}(\tau) 
\end{equation*}
Thus the  Cayley graph $Hyp_f(1-t)$ is exactly  isomorphic to the Cayley
graph of $\psi(\Gamma_{1-t}(\tau))$ with respect to the following generators:
\begin{equation}\label{mat888}
 \psi(\prod_{j=0}^i (1 - x_{d-j}^{1-q}\tau)) =\prod_{j=1}^i \psi(1 - x_{d-j}^{1-q}\tau) \end{equation}
$ \text{where for any }\quad i=1,\ldots, d-1 $
\begin{equation*}
 \psi(\varphi_j)=\begin{pmatrix}  
1 &-x_{d-i}^{1-q}  &  & \ldots & 0\\
                 0 & 1 & -x_{d-i}^{q-q^2} & \ldots &
                \\\vdots & \vdots &
                 \vdots & \ddots & \vdots \\ 0 & 0 & 0 & \ldots &-x_{d-i}^{q^{d-2}-q^{d-1}} \\-tx_{d-i}^{q^{d-1}-q^d} &
                 0 & 0 & \ldots & 1 \ 
\end{pmatrix},
\end{equation*}
for all elements  $(x_1,\ldots,x_{d-1}) \in \mathcal{B}_{q,d}^{1-t}$.
\begin{definition}
Motivated by the  classical Legendre symbol, we define here an extended  Legendre
symbols for arbitrary $d$ and $A, B \in \mathbb{F}_q[t]$ $($or in $\mathbb{Z})$ by:\begin{equation}\label{Leg}
\biggl(\frac{A}{B}\biggr)_d : = \begin{cases}1 & \quad \text{if}\quad {X^d \equiv A
    \mod{B}}\quad \text{has a solution}\\
                                           -1 & \quad otherwise
\end{cases}
\end{equation}
\end{definition}
\begin{theorem}
Let q be a an odd prime power and $f(t)$ be  irreducible of degree equal
to $dn$ for some positive integer $n$, and let $\psi$ be as $($\ref{pepsi}$)$. Then:
\begin{equation*}
\text{Image}(\psi) = \begin{cases} 
\text{PSL}(d, \mathbb{F}_{q^{\deg{f}}}) &
    \qquad  \text{if}\quad \biggl(\frac{1-t}{f(t)}\biggr)_d = 1\\
   \text{PGL}(d, \mathbb{F}_{q^d}) & \qquad \text{otherwise}
\end{cases}
\end{equation*}
\end{theorem} 
\begin{proof} 
Set 
\[
 \mathcal{U} : = D^{(1)}(\mathbb{F}_q[t]_\infty)\prod_{h\neq t, 1-t}J_h
\Gamma'^{(1)}(\tau)
\]

where  $D^{(1)}$ given by \ref{D1} and 
\[
 J_f : =   \ker{\biggl(D^{(1)}(\mathbb{F}_q[t]_h)\quad \longrightarrow \quad D^{(1)}(\mathbb{F}_q[t]_h/f\mathbb{F}_q[t]_f)\biggr)}
\]
and 

\[
 J_h : =   D^{(1)}(\mathbb{F}_q[t]_h)\qquad \text{if} \quad h \neq f
\]
and finally 
\[\Gamma'^{(1)}(\tau):=\ker{\biggl(\mathbb{F}_{q^d}\{\tau\}[\frac{1}{1-t}]^\times\quad
  \longrightarrow \quad \mathbb{F}_q\{\{\tau\}\}^\times}\biggr).
\]
 Define 
\[\widetilde{\psi} : \mathcal{U}\quad \longrightarrow \quad
\text{PGL}(d,\mathbb{F}_{q^{\deg{f}}}) \]
as composition $\text{Proj}_f\text{mod}_f$ where
\[( \zeta_\infty,\ldots,\zeta_f,\ldots)\quad \mapsto \quad \zeta_f,\]
By Theorem \ref{96} we know that $D = \mathbb{F}_{q^d}(\tau)$ is
unramified at $f$. Thus we have an isomorphism   $D_f
=D\otimes_{\mathbb{F}_q(t)}\mathbb{F}_q(t)_f \cong \text{M}(d,\mathbb{F}_q(t)_f)$
which takes the maximal order $\mathbb{F}_{q^d}\{\tau\}\otimes_{\mathbb{F}_q[t]}\mathbb{F}_q[t]_f$ to
the maximal order $\text{M}(d,\mathbb{F}_q[t]_f)$. So
\[\mathbb{F}_{q^d}\{\tau\}\otimes_{\mathbb{F}_q[t]}\mathbb{F}_q[t]_f/f\mathbb{F}_q[t]_f \cong \text{M}(d, \mathbb{F}_q[t]_f/f\mathbb{F}_q[t]_f)\]
and this induces an isomorphism : 
\begin{align*}
(\mathbb{F}_{q^d}\{\tau\}\otimes_{\mathbb{F}_q[t]}\mathbb{F}_q[t]_f/f\mathbb{F}_q[t]_f)^{\times}/Z
\quad &\cong \quad \text{PGL}(d,
\mathbb{F}_q[t]_f/f\mathbb{F}_q[t]_f)\\
 \cong &\quad \text{PGL}(d, \mathbb{F}_{q^{\deg{f}}})
\end{align*}
Thus $\widetilde{\psi}$ acts as projection on the $f^{\text{th}}$ component and after
reduction modulo $f$ sends this component, under the above isomorphism $($which can be exactly
our $\psi$$)$ to the element  defined by $($\ref{mat888}$)$. Again applying Theorem \ref{96} we see
that 
\[ D^{(1)}(\mathbb{F}_q[t]_{1-t}) \cong \text{SL}(d,  \mathbb{F}_q[t]_{1-t}),\]

and so from Strong approximation Theorem,  it follows
immediately that \[ D^{(1)}(\mathbb{F}_q[t]) D^{(1)}(\mathbb{F}_q[t]_{1-t}) \] is
dense in $ D^{(1)}(\mathbbb{A})$, i.e, for the  open set $\mathcal{U}$ $($and for all open sets$)$ of $D^{(1)}(\mathbbb{A}) /
D^{(1)}(\mathbb{F}_q[t]_{1-t})$, 
$D^{(1)}(\mathbb{F}_q[t]) \cap \mathcal{U}$ must be dense in $\mathcal{U}$.\\ 
Now since $\widetilde{\psi}$ is a continuous function over $\mathcal{U}$ with a
finite range, we must have :
\[ \widetilde{\psi}(D^{(1)}(\mathbb{F}_q[t]) \cap \mathcal{U}) =
 \widetilde{\psi}(\mathcal{U}).\]
 We can see directly that
\[ D^{(1)}(\mathbb{F}_q[t]) \cap \mathcal{U} =
\{ \mu \in \Gamma'^{(1)}(\tau)\quad | \quad  \mu \equiv 1 \mod{f}\}.\]
So this shows that  
\begin{align*}
\widetilde{\psi}(D^{(1)}(\Gamma^{1-t}(\tau))) &=
\psi(D^{(1)}(\Gamma^{1-t}(\tau))) \supseteq  \widetilde{\psi}(\{ \mu \in
\Gamma'^{(1)}(\tau)\quad | \quad  \mu \equiv 1 \mod{f}\} \\&= \text{PSL}(d,\mathbb{F}_{q^{\deg{f}}}) 
\widetilde{\psi}(\mathcal{U}).
\end{align*}
 since 
\[X \in \text{PSL}(d,\mathbb{F}_{q^{\deg{f}}}) \quad \text{iff}
 \quad \biggl(\frac{\det{X}}{f(t)}\biggr)_d = 1.\]
It is an immediate
 result of the following diagram:
\begin{diagram}
&& 1 && 1 && 1 && \\
&& \dTo && \dTo && \dTo && \\
1 & \rTo & U_d & \rTo & k^\times & \rTo^{ x \mapsto x^d} &  k^{\times d} & \rTo & 1 \\
&& \dTo && \dTo && \dTo && \\
1 & \rTo & \text{SL}(d; k) & \rTo & \text{GL}(d; k) & \rTo^\det &  k^\times  & \rTo
& 1 \\
&& \dTo && \dTo && \dTo && \\
1 & \rTo & \text{PSL}(d; k) & \rTo & \text{PGL}(d; k) & \rTo &  k^\times /
k^{\times d} & \rTo & 1 \\
&& \dTo && \dTo && \dTo && \\
&& 1 && 1 && 1 && \\
\end{diagram}
where $ k = \mathbb{F}_{q^{\deg{f}}}$ and $ U_d := \{ u \in k \quad | \quad u^d =
1\}$. Exactness of columns and rows and commutativity of the diagram $($for
arbitrary k$)$ is well known.
We see by assumption that any generator of $\Gamma^{1-t}(\tau)$ is in $
\text{PSL}(d,\mathbb{F}_{q^{\deg{f}}})$ iff $\biggl(\frac{1-t}{f(t)}\biggr)_d =
1$. The other case will be handled exactly as the known case for $d = 2$, see
  \cite[Theorem $4.13$] {mor1}. 
\end{proof}
\begin{cor}
Associated to  any element  $(x_1,\ldots,x_{d-1}) \in
\mathcal{B}_{q,d}^{1-t}$  are the matrices:
\begin{equation}\label{mat999}
 \text{M}_{i,\ldots,1} :=\prod_{j=1}^{d-1} \text{Q}_j \qquad
\end{equation}
 $\text{where for any }\quad j=1,\ldots, d-1$ 
 \begin{equation*}
 \text{Q}_j=\begin{pmatrix}  
1 &-x_{d-j}^{1-q}  &  & \ldots & 0\\
                 0 & 1 & -x_{d-j}^{q-q^2} & \ldots &
                \\\vdots & \vdots &
                 \vdots & \ddots & \vdots \\ 0 & 0 & 0 & \ldots &-x_{d-j}^{q^{d-2}-q^{d-1}} \\-tx_{d-j}^{q^{d-1}-q^d} &
                 0 & 0 & \ldots & 1
\end{pmatrix}
\end{equation*}
\begin{equation*}
 i) \quad \text{If}\quad  \biggl(\frac{1-t}{f(t)}\biggr)_d = 1\end{equation*}
Then $Hyp_f(1-t)$ is exactly  isomorphic to the Cayley
graph of $\text{PSL}(d, \mathbb{F}_{q^{\deg{f}}})$ with respect to the
generators \ref{mat999}.\\
\begin{equation*}
ii) \quad \text{If}\quad  \biggl(\frac{1-t}{f(t)}\biggr)_d = -1\end{equation*}
Then $Hyp_f(1-t)$ is exactly  isomorphic to the Cayley
graph of $\text{PGL}(d, \mathbb{F}_{q^{\deg{f}}})$ with respect to the
generators \ref{mat999}.\end{cor}
\begin{proof}
 This follows immediately from the above Theorem.
\end{proof}